\documentclass[a4paper,11pt]{article}

\usepackage[textwidth=6.4in,textheight=9.0in]{geometry}

\usepackage{amsthm}
\usepackage{amsmath}
\usepackage{amssymb}
\usepackage{graphicx}
\usepackage[round]{natbib}
\usepackage{bbm}
\usepackage{textcomp}
\usepackage[all]{xy}

\usepackage{color}
\usepackage{framed}
\usepackage{comment}
\definecolor{shadecolor}{gray}{0.9}

\specialcomment{extra}{\begin{shaded}}{\end{shaded}}

\newlength{\myleftmargin}
\usepackage{amsmath}
\usepackage{amssymb}
\usepackage{amsthm}
\usepackage{color}
\DeclareSymbolFontAlphabet{\Bbb}{AMSb}

\newtheorem{theorem}{Theorem}[section]
\newtheorem{lemma}[theorem]{Lemma}
\newtheorem{proposition}[theorem]{Proposition}
\newtheorem{corollary}[theorem]{Corollary}
\newtheorem{definition}[theorem]{Definition}

\newtheorem{remark}[theorem]{Remark}

\newcommand{\atob}[2]{\emph{#1)} $\Rightarrow$ \emph{#2)}.} 
\newcommand{\aeqb}[2]{\emph{#1)} $\Leftrightarrow$ \emph{#2)}.} 
\newcommand{\ada}[1]{\emph{#1).}}

\newenvironment{proofof}[1]{\noindent{\bf Proof of #1:}}{\qed\medskip}

\newlength{\fixboxwidth}
\definecolor{darkgreen}{rgb}{0,0.6,0}

\newcommand{\ca}[1]{{\cal #1}}

\newcommand{\C}{\mathbb{C}}
\newcommand{\N}{\mathbb{N}}

\newcommand{\R}{\mathbb{R}}
\newcommand{\T}{\mathbb{T}}
\newcommand{\Rd}{\mathbb{R}^d}

\newcommand{\Z}{\mathbb{Z}}

\newcommand{\E}{\mathbb{E}}

\DeclareMathOperator{\re}{Re}
\DeclareMathOperator{\im}{Im}
\newcommand{\iu}{\mathrm i}

\renewcommand{\a}{\alpha}

\newcommand{\g}{\gamma}

\renewcommand{\d}{\delta}
\newcommand{\D}{\Delta}
\newcommand{\e}{\varepsilon}

\renewcommand{\k}{\kappa}
\newcommand{\lb}{\lambda}

\newcommand{\s}{\sigma}

\renewcommand{\t}{\tau}
\newcommand{\p}{\varphi}
\renewcommand{\P}{\Phi}
\newcommand{\om}{\omega}

\DeclareMathOperator{\codim}{codim}

\DeclareMathOperator{\spann}{span}

\DeclareMathOperator{\ran}{ran}

\DeclareMathOperator{\id}{id}
\newcommand{\eins}{\boldsymbol{1}}

\newcommand{\snorm}[1]{\Vert #1 \Vert}
\newcommand{\mnorm}[1]{\bigl\Vert \, #1 \, \bigr\Vert}

\newcommand{\hnorm}[1]{\Vert #1 \Vert_{H}}
\newcommand{\Hnorm}[1]{\Bigl\Vert #1 \Bigr\Vert_{H}}
\newcommand{\inorm}[1]{\Vert #1 \Vert_\infty}
\newcommand{\tvnorm}[1]{\Vert #1 \Vert_{\mathrm TV}}

\newcommand{\sLx}[2]{{\ca L_{#1}(#2)}}

\newcommand{\sLtn}{\sLx 2 \nu}

\newcommand{\Lx}[2]{{L_{#1}(#2)}}

\newcommand{\Ltn}{\Lx 2 \nu}

\newcommand{\Tkx}[1]{T_{k,#1}}
\newcommand{\Ikx}[1]{I_{k,#1}}
\newcommand{\Ikxs}[1]{I_{k,#1}^*}

\newcommand{\Tkn}{\Tkx \nu}
\newcommand{\Ikn}{\Ikx \nu}
\newcommand{\Ikns}{\Ikxs \nu}

\newcommand{\sosbmh}[3]{\ca M_{#1}^{#2}(#3)}
\newcommand{\sosbm}[1]{\sosbmh {}{}{#1}}
\newcommand{\sobm}[1]{\sosbmh {+}{}{#1}}
\newcommand{\sobpm}[1]{\sosbmh {1}{}{#1}}

\newcommand{\sosrmh}[3]{\ca M_{#1}^{#2}(#3)}
\newcommand{\sosrm}[1]{\sosrmh {}{*}{#1}}
\newcommand{\sorm}[1]{\sosrmh {+}{*}{#1}}

\newcommand{\aec}[1]{[{#1}]_\sim}

\theoremstyle{plain}  
\newtheorem{thm}{Theorem}[section]

\theoremstyle{definition} 
\newtheorem{defn}{Definition}[section]

\theoremstyle{remark}

\newcommand{\diff}{\,\mathrm{d}}

\newcommand{\bbS}{\mathbb{bbS}}

\defcitealias{dlmf}{DLMF}

\begin{document}

\title{Strictly proper kernel scores and characteristic kernels on compact spaces}
\author{Ingo Steinwart \and Johanna F.~Ziegel}

\maketitle

\begin{abstract}
Strictly proper kernel scores are well-known tool in probabilistic forecasting, while characteristic
kernels have been extensively investigated in the machine learning literature. We first show that 
both notions coincide, so that insights from one part of the literature can be used in the other.
We then show that the metric induced by a characteristic kernel cannot reliably distinguish 
between distributions that are far apart in the total variation norm as soon as the underlying 
space of measures is infinite dimensional. In addition, we provide a characterization 
of characteristic kernels in terms of eigenvalues and -functions and apply this 
characterization to the case of continuous kernels on (locally) compact   spaces.
In the compact case we further show that characteristic kernels exist if and only if 
the space is metrizable.
As special cases of our general theory we investigate translation-invariant kernels on compact Abelian groups
and isotropic kernels on spheres. The latter are of particular interest for forecast evaluation of probabilistic predictions on spherical domains as frequently encountered in meteorology and climatology.
\end{abstract}

\section{Introduction}\label{sec:intro}

Probabilistic forecasts 
of uncertain future events
are issued in a wealth of applications, see \citet{GneitingKatzfuss2014} and the references therein. To 
assess the quality and to compare such forecasts,
proper scoring rules are a well-established tool, see \citet{GneitingRaftery2007},
and in applications, it is usually even 
desirable to work with  \emph{strictly} proper scoring rules.
 A broad class of proper scoring rules are so-called kernel scores, which are constructed using a positive definite kernel. Unfortunately, 
 however, no general conditions are available to decide whether a given kernel induces a strictly proper kernel score. As detailed below in Theorem \ref{thm:one}, strict propriety of a kernel score is intimately connected to the kernel being characteristic, a notion that has been  studied in the machine learning literature for a decade,
see e.g.~\cite{GrBoRaScSm07a, FuSrGrSc09a, SriperumbGrettonETAL2010, SriperumbFukumizuETAL2011} as well as 
the recent survey of \cite{MuFuSrSc17a} and the references therein.
 In this paper, we study characteristic kernels on compact spaces extending results of \citet{MiXuZh06a} and \citet{SriperumbGrettonETAL2010,SriperumbFukumizuETAL2011}. As a consequence, we can characterize strictly proper kernel scores on compact Abelian groups and the practically highly relevant example of spheres.

To describe our results in more detail, let us formally introduce some of the notions mentioned above.
To this end, 
let $(X,\mathcal{A})$ be a measurable space and let $\mathcal{M}_1(X)$ denote the class of all probability measures on $X$. For $\mathcal{P}\subseteq \mathcal{M}_1(X)$, a \emph{scoring rule} is a function $S:\mathcal{P}\times X \to [-\infty,\infty]$ such that the integral $\int S(P,x) \diff Q(x)$ exists for all $P, Q \in \mathcal{P}$. The scoring rule is called \emph{proper} with respect to $\mathcal{P}$ if
\begin{equation}\label{eq:proper}
\int S(P,x) \diff P(x) \le \int S(Q,x) \diff P(x), \quad \text{for all $P,Q \in \mathcal{P}$},
\end{equation}
and it is called \emph{strictly proper} if equality in \eqref{eq:proper} implies $P=Q$. 
Recall that if
 the class $\mathcal{P}$ consists of absolutely continuous probability measures with respect to some $\sigma$-finite measure $\mu$ on $X$ then the logarithmic score 
$S(P,x) := -\log p(x)$, where $p$ is the density of $P$, is a widely 
used example of a strictly proper scoring rule for density forecasts. 
Another well-known example is the  Brier score for distributions on 
$X = \{1,\dots,m\}$ that is 
 defined as $S(P,i) := \sum_{j=1}^m p_j^2 + 1 - 2p_i$, where $p_i = P(\{i\})$, $i=1,\dots,m$. Finally, for $X = \R$, the continuous ranked probability score (CRPS) is given by 
\[
S(P,x) := \int_\R |y - x|\diff P(y) - \frac{1}{2} \int_\R \int_\R |y - y'|\diff P(y)\diff P(y').
\]
It is strictly proper with respect to the class of all probability measures 
with finite first moment, see e.g.~\citet[Section 4.2]{GneitingRaftery2007}, and 
consequently 
 it can be used to evaluate predictions of density 
forecasts as well as probabilistic forecasts of categorical variables. Various other examples 
can be found in \citet{GneitingRaftery2007}. 

One general class of proper scoring rules are kernel scores. 
To this end, let 
$k:X \times X \to \mathbb{R}$ be a symmetric function. We call $k$ a \emph{kernel}, if 
 it is positive definite, that is, if 
\begin{equation}\label{posdef}
   \sum_{i,j=1}^n a_i a_j k(x_i,x_j) \ge 0
\end{equation}
 for all natural numbers $n$, all $a_1,\dots,a_n \in \R$, and all $x_1,\dots,x_n \in X$. 
It is strictly positive definite if equality in \eqref{posdef} 
implies $a_1=\dots=a_n=0$ whenever the points  $x_1,\dots,x_n$ are mutually distinct.
Let us assume that $k$ is measurable and define
\begin{displaymath}
   \mathcal{M}_1^k(X):=\Bigl\{P \in \mathcal{M}_1(X)\;\bigl|\; \int_{X}\sqrt{k(x,x)} \diff P(x) < \infty\Bigr\}.
\end{displaymath}
For a bounded kernel $k$, we have that $\mathcal{M}_1^k(X) = \mathcal{M}_1(X)$,
and the Cauchy-Schwarz inequality
 $k(x,y) \le \sqrt{k(x,x)}\sqrt{k(y,y)}$ for kernels 
shows that, for all $P,Q \in \mathcal{M}_1^k(X)$, the kernel $k$ is integrable with respect to the product measure $P\otimes Q$. 

\begin{defn}
The \emph{kernel score} $S_k$ associated with a measurable kernel $k$ on $X$ 
is the scoring rule $S_k:\mathcal{M}_1^k(X)\times X \to \mathbb{R}$ defined by
\begin{equation*}
S_k(P,x) := - \int k(\omega,x) \diff P(\omega) + \frac{1}{2}\int \int k(\omega,\omega')\diff P(\omega)\diff P(\omega').
\end{equation*} 
\end{defn}

Kernel scores are a broad generalization of the CRPS, and  in fact, also the Brier score
 can be rewritten as a kernel score, see \citet[Section 5.1]{GneitingRaftery2007}.
 However, the logarithmic score does not belong to this class. 
If
 $X$ is a Hausdorff space and $k$ is continuous, 
then \citet[Theorem 4]{GneitingRaftery2007} show that
 $S_k$ is proper with respect to all Radon probability measures on $X$. 
Their result is based on \citet[Theorem 2.1, p.~235]{BergChristensETAL1984}, 
where it is fundamental that the kernel is continuous. In this respect 
we remark that the definition of a kernel score of \citet{GneitingRaftery2007} 
is more general than ours as it allows for kernels being only conditionally positive definite. 
While this level of generality is fruitful for example in the case $X = \R^d$, 
we believe that it is sufficient to consider only positive definite kernels 
for compact spaces. Indeed, note that if $X$ is compact and separable and there is a
 strictly positive probability measure $\nu$ on its Borel sets, then, 
if we only consider kernels $k$ such that $\int_X k(x,y)\diff \nu(x)$ does
 not depend on $y$, we know by \citet[Theorem 2]{Bochner1941}
 that conditionally positive definite kernels and 
positive definite kernels are the same up to a constant.

In our framework, it is possible to show propriety of the kernel score without requiring continuity of the kernel using the theory of reproducing kernel Hilbert spaces (RKHS); see Theorem \ref{thm:one} below. In addition, we obtain a condition for when the kernel score is actually strictly proper. 

\begin{thm}\label{thm:one}
Let $k$ be a measurable kernel with RKHS $H$ with norm $\|\cdot\|_H$ and 
$\Phi: \mathcal{M}_1^k(X) \to H$ be the kernel embedding defined by 
\begin{equation}\label{eq:inj}
\Phi(P) :=  \int k(\cdot,\omega)\diff P(\omega)\, .
\end{equation}
Then the kernel score satisfies 
\begin{equation}\label{main-obs}
   \mnorm{\P(P)-\P(Q)}_H^2 = 2 \left(\int S_k(Q,x)\diff P(x) - \int S_k(P,x)\diff P(x)\right)
\end{equation}
for all $P,Q\in \mathcal{M}_1^k(X)$. In particular, 
 $S_k$ is a proper scoring rule with respect to $\mathcal{M}_1^k(X)$, and it  is strictly proper if and only if 
$\P$ is injective.
\end{thm}

In the machine learning literature a bounded
measurable kernel is called 
 \emph{characteristic} if the kernel (mean) embedding $\Phi: \mathcal{M}_1(X) \to H$
defined by \eqref{eq:inj} is injective.
Consequently, 
Theorem \ref{thm:one} shows that, for bounded measurable kernels $k$, strictly proper kernel
scores $S_k$ are exactly those for which $k$ is characteristic.
In particular, the wealth of examples and conditions for characteristic kernels
can be directly used to find new strictly proper scoring rules and vice versa.

While Theorem \ref{thm:one} is an interesting observation for 
both machine learning applications and  probabilistic forecasting, its proof is actually 
rather trivial. In the rest of the paper, we therefore focus on more involved aspects of characteristic kernels,
and strictly proper kernel scores, respectively.
We introduce the necessary mathematical machinery 
on kernels and their interaction with (signed) finite  measures 
in Section \ref{sec:prelim}.
In particular we recall that 
 a bounded measurable kernel $k$ with RKHS $H$ induces a semi-norm $\|\cdot\|_H$ on $\mathcal{M}(X)$, the space of all finite signed measures on $X$, via the kernel mean embedding, that is via the left-hand side of \eqref{main-obs}.
 In Section \ref{sec:general}, we study this semi-norm for general $X$. 
In particular, Theorem \ref{norm-equivalence} shows that for injective kernel embeddings,
$\|\cdot\|_H$ fails to be equivalent to the 
total variation norm
on $\mathcal{M}(X)$ if and only if $\dim \mathcal{M}(X)=\infty$, and Corollary
\ref{different-distances-cor} gives an even sharper result on $\mathcal{M}_1(X)$.
In view of \eqref{main-obs}, these results show that the value  on the right-hand side of
\eqref{main-obs} is not proportional to the (squared) total variation norm. 
Besides some structural results on characteristic kernels, see Lemmas \ref{char-add} and    \ref{product-form},
we further present a simple computation of the left-hand side of \eqref{main-obs}
in terms of eigenvalues and -functions of a suitable integral operator.
In Section \ref{sec:lcs}, we then  exploit our general theory to obtain new results for bounded continuous kernels on locally compact Hausdorff spaces. The main question of interest is when such kernels are universal or characteristic. 
In Theorem \ref{thm-new-char-3} and Corollary \ref{uni-cor} we
give characterization results in terms of eigenfunctions of certain integral operators.
We also provide insight concerning the difference of considering Borel or Radon measures on locally compact Hausdorff spaces in the study of kernel embeddings, see Theorems \ref{exist-sipd} and \ref{thm-new-char-2}. As a result, it turns out in Theorem \ref{exists-univer-char-kern} that continuous characteristic kernels on compact Hausdorff spaces only exist if the spaces are metrizable. In Section \ref{sec:structure}, we apply the characterization results of Section \ref{sec:lcs} to 
translation-invariant kernels on 
compact Abelian groups and to isotropic kernels on spheres. 
All proofs can be found in Section \ref{sec:proofs}.



\section{Preliminaries}\label{sec:prelim}

In this section we recall some facts about reproducing kernels and their interaction with measures.
To this end, let $(X,\ca A)$ be a measurable space. We denote the space of 
finite signed measures on $X$ by $\sosbm X$  and write 
$\sobm X$ and $\sobpm X$ for the subsets of all (non-negative) finite, respectively probability measures.
Moreover, we write 
\begin{displaymath}
 \sosbmh 0 {} X := \bigl\{ \mu \in \sosbm X: \mu(X) = 0\bigr\}\, .
\end{displaymath}
As usual we equip $\sosbm X$ and its subsets above
with the total variation norm $\tvnorm\cdot$. 
Recall that $\tvnorm\cdot$ is complete and hence  $\sosbm X$ is a Banach space.
Moreover, $\sosbmh 0 {} X$ is closed subspace of co-dimension 1, which contains, e.g.~all differences
of probability measures. Moreover, for every $P\in \sosbmh 1 {} X$ we have 
\begin{equation}\label{m0+p}
 \sosbm X = \R P \oplus \sosbmh 0 {} X\, .
\end{equation}

Given a measurable function $f:X\to \R$ and a measure $\nu$ on $X$ we write 
$[f]_\sim$ for the $\nu$-equivalence class of $f$. Similar, we denote the space of 
$p$-times $\nu$-integrable functions by $\sLx p \nu$ and the corresponding 
space of $\nu$-equivalence classes by $\Lx p \nu$. Note that this rather pedantic
notation becomes very useful when dealing with RKHSs since these spaces consist
of functions that can be  evaluated pointwise in a continuous fashion, whereas 
such an evaluation does, in general, not make sense for elements in $\Lx p \nu$.

To formally introduce kernel mean embeddings, we need to recall the notion of 
Pettis integrals. To this end let, let $H$ be a Hilbert space and $f:X\to H$ be a function.
Then $f$ is weakly measurable, if $\langle w, f\rangle :X\to \R$ is measurable for all $w\in H$.
Similarly, $f$ is weakly integrable with respect to a measure $\nu$ on $(X,\ca A)$,
if $\langle w, f\rangle\in \sLx 1 \nu$ for all $w\in H$. In this case, 
there exists a unique
$i_\nu(f)\in H$, called the Pettis integral of $f$ with respect to $\nu$, such that for all $w\in H$ we have 
\begin{displaymath}
   \langle w, i_\nu(f) \rangle = \int_X \langle w, f\rangle\, \diff\nu\, ,
\end{displaymath}
see e.g.~\citet[Chapter II.3]{DiUh77} together with the reflexivity of $H$ and the identity $H=H'$
between $H$ and its dual $H'$. Using the Hahn-Jordan decomposition, it is not hard to see that 
$i_\mu(f)$ can analogously be defined for finite signed measures $\mu$. In the following, we adopt the more intuitive
notation $\int_X f \diff\mu := i_\mu(f)$, so that the defining equation above  becomes
\begin{align}\label{pettis-formula}
    \bigl\langle w, \int_X f \diff\mu\bigr\rangle = \int_X \langle w, f\rangle\, \diff\mu\, .
\end{align}
Furthermore, in the case of probability measures $\mu$, we sometimes also write $\E_\mu f := i_\mu(f)$.
Let us now use Pettis integrals to define kernel mean embeddings. To this end, let $H$
be an RKHS over $X$ with kernel $k$ and canonical feature map $\P:X\to H$, that is 
$\P(x) := k(\cdot,x)$ for all $x\in X$. Then $\P$ is weakly measurable if and only if 
$\langle h , \P\rangle = h$ is measurable for all $h\in H$, and therefore, we conclude that $\P$
is weakly integrable with respect to some measure $\nu$ on $X$, 
if and only if $h\in \sLx 1 \nu$ for all $h\in H$.
By a simple application of the closed graph theorem, the latter is equivalent to the 
continuity of the map $[\, \cdot\, ]_\sim :H\to \Lx 1 \nu$.
In this respect recall that $H$ consists of measurable functions if and only if $k$ is separately measurable, 
that is
$k(\cdot, x):X\to \R$ is measurable for all $x\in X$, see \citet[Lemma 4.24]{StCh08},
and $[\, \cdot\, ]_\sim :H\to \Lx 1 \nu$ is continuous if, e.g.
\begin{align}\label{square-finite}
   \int \sqrt{k(x,x)} \diff\nu(x)<\infty\, ,
\end{align}
see e.g.~\citet[Theorem 4.26]{StCh08}. Obviously, the latter condition is still sufficient for 
finite signed measures if one replaces $\nu$ by $|\nu|$.
For a separately measurable kernel $k$ with RKHS $H$ we now write  
\begin{displaymath}
   \sosbmh {} k X := \Bigl\{ \mu\in \sosbm X: H\subset \sLx 1\mu \Bigr\} \, ,
\end{displaymath}
and analogously we define $\sosbmh {+} k X $ and $\sosbmh {1} k X $. Obviously, $\sosbmh {} k X$ is 
a vector space containing all Dirac measures and 
using \eqref{square-finite} it is not hard to see that $\sosbmh {} k X = \sosbmh {} {} X$,
if $k$ is bounded. In fact, the latter is also necessary for $\sosbmh {} k X = \sosbmh {} {} X$ as a combination
 of  \citet[Proposition 2]{SriperumbGrettonETAL2010} with \citet[Lemma 4.23]{StCh08} shows.
 Moreover, our considerations 
above show that  $\sosbmh {} k X $ is the largest set on which we can define the kernel embedding 
\begin{displaymath}
   \P(\mu) := \int_X \P \diff\mu\ = \int_X k(\cdot, x)\diff \mu(x)\, .
\end{displaymath}
Note that the map $\P: \sosbmh {} k X \to H$
is linear, and 
for Dirac measures $\d_x$ we have $\P(\d_x) = \P(x)$.
Consequently  
\begin{displaymath}
   \hnorm \mu := \hnorm { \P(\mu) }
\end{displaymath}
defines a new semi-norm on $\sosbmh {} k X$, and by a double application of 
\eqref{pettis-formula} we further have 
\begin{align}\label{hinner}
   \langle \P(\mu_1) , \P(\mu_2)\rangle_H 
=  \int_X \int_X k(x,x') \diff\mu_1(x)\diff \mu_2(x') \, .
\end{align}
Note that this semi-norm is a norm, if and only if 
the kernel embedding $\P: \sosbmh {} k X \to H$ is injective
and by \eqref{hinner} the latter is equivalent to 
\begin{align}\label{sipd}
 \snorm\mu_H^2 = \int \int k(x,x') \diff \mu(x)\diff \mu(x') >0 \, .
\end{align}
for all $\mu\in \sosbmh {} k X\setminus \{0\}$. 
This leads to the following definition.

\begin{definition}
  Let $k$ be a  measurable kernel on $X$ and 
$\ca M\subset 
\sosbmh {} k X$. Then $k$ is 
 called strictly
integrally positive definite with respect to   $\ca M$, if \eqref{sipd}
holds for all $\mu\in \ca M$ with $\mu\neq 0$.
\end{definition}

It is well-known, that using \eqref{pettis-formula}
 the  semi-norm introduced above can also  be computed by 
\begin{align}\nonumber
   \hnorm \mu 
= \Hnorm { \int_X \P \diff\mu } 
= \sup_{f\in B_H} \Bigl|\bigl\langle f , \int_X \P \diff\mu\bigr\rangle   \Bigr|  
&= \sup_{f\in B_H} \Bigl| \int_X \langle f , \P \rangle \diff\mu   \Bigr|  \\  \label{alternative-comp}
&= \sup_{f\in B_H} \Bigl| \int_X f \diff\mu   \Bigr|  \, , 
\end{align}
where $B_H$ denotes the closed unit ball of $H$. 
Consequently, we have  
$\hnorm \mu \leq \snorm{[\, \cdot\, ]_\sim :H\to \Lx 1 \mu}$,
and if $k$ is bounded, then 
$\P: \sosbmh {} {} X \to H$ is continuous with 
$\snorm {\P: \sosbmh {} {} X \to H} \leq \inorm{k}^{1/2}$. 
In particular, if  $k$ is bounded and $\P: \sosbmh {} {} X \to H$ is injective,
then $\hnorm\cdot$ defines a new norm on $\sosbm X$ that is dominated by 
$\tvnorm \cdot$ and that describes a Euclidean geometry with inner product \eqref{hinner}.
Unless $\dim  \sosbm X< \infty$, however, both norms are not equivalent, as we will see in Theorem \ref{norm-equivalence}.

With the help of the new semi-norm $\hnorm \cdot$ on  $\sosbmh {} k X$
we can now define a semi-metric on $\sosbmh 1 k X$ by setting 
\begin{displaymath}
   \g_k(P,Q) :=  \hnorm{P-Q} 
= \sup_{f\in B_H} \left|\int f\diff P - \int f\diff Q\right| 
\end{displaymath}
for $P,Q\in \sosbmh 1 k X$. Here, we note that the second equality, which 
follows from our considerations above, has already been 
shown by \citet[Theorem 1]{SriperumbGrettonETAL2010}. 
Similarly, the following definition is taken from 
\citet{FuGrSuSc08a, SriperumbGrettonETAL2010}.

\begin{definition}\label{def-char-kern}
   A bounded measurable kernel $k$ on $X$ is called characteristic, if 
the kernel mean embedding 
$\P_{|\sosbmh 1 {} X}:\sosbmh 1 {} X\to H$ is injective.
\end{definition}

%
Clearly, $k$ is characteristic, if 
and only if $\g_k$ is a metric, and a literal repetition of \citet[Lemma 8]{SriperumbGrettonETAL2010}
shows:

\begin{proposition}\label{M0-char}
 Let $(X,\ca A)$ be a measurable space and $k$ be a bounded measurable kernel on $X$. Then the following statements
 are equivalent:
 \begin{enumerate}
 \item $k$ is strictly integrally positive definite with respect to $\sosbmh 0 {} X$.
  \item $k$ is characteristic.
 \end{enumerate}
\end{proposition}


%

Now, let $\nu$ be a measure  on $X$, $k$ be a measurable kernel on $X$,
and $H$ be its RKHS.  We further assume that the map $\Ikn:H\to \Lx 2 \nu$ 
given by $f\mapsto [f]_\sim$  is well-defined and compact. For an example of 
such a situation recall \citet[Lemma 2.3]{StSc12a}, which shows that $\Ikn$ is Hilbert-Schmidt if
\begin{align}\label{int-diag}
   \int_X k(x,x) \, \diff\nu(x) < \infty\, .
\end{align}
Obviously, the latter holds, if e.g.~$k$ is bounded and $\nu$ is a finite measure.
Now assume that $\Ikn$ is well-defined and compact. Then, the 
associated integral operator  $\Tkn:\Ltn \to \Ltn$, 
defined by 
\begin{equation*}
   \Tkn f = \Bigl[\int_X k(x, \cdot) f(x)\, \diff\nu(x)\Bigr]_\sim\, , \qquad f\in \Ltn\, ,
\end{equation*}
satisfies $\Tkn = \Ikn \circ \Ikns$, see e.g.~\citet[Lemma 2.2]{StSc12a}, where 
$\Ikns$ denotes the adjoint of $\Ikn$. In particular, 
$\Tkn$ is compact, 
positive, and self-adjoint, and if \eqref{int-diag} is satisfied, $\Tkn$ is even nuclear.
%
%
Moreover,
 the spectral theorem in the form of \citet[Lemma 2.12]{StSc12a} 
gives us an at most countable, ordered family  $(\lb_i)_{i\in I}\subset (0,\infty)$ converging to $0$
and a family $(e_i)_{i\in I}\subset H$ such that:
\begin{itemize}
   \item $(\lb_i)_{i\in I}$ are the non-zero eigenvalues 
of $\Tkn$ including multiplicities,
	\item  $([e_i]_\sim)_{i\in I}$
is an $\Ltn$-ONS of the corresponding eigenfunctions with 
\begin{align}\label{span_ei}
   \overline{\spann\{[e_i]_\sim : i\in I \}}^{\Ltn} = \overline{ [H_\sim]}^{\Ltn}\, ,
\end{align}
	\item $(\sqrt{\lb_i} e_i)_{i\in I}$ is an ONS in $H$.
\end{itemize}
Here, 
we say that an at most countable family $(\a_i)_{i\in I} \subset [0,\infty)$
converges to 0, if either $I=\{1,\dots,n\}$ or $I=\N$ and $\lim_{i\to \infty}\a_i = 0$.

Note that for nuclear $\Tkn$, we additionally have $\sum_{i\in I}\lb_i < \infty$, and if 
$k$ is bounded, $\inorm {e_i}<\infty$ holds for all $i\in I$.
Finally, \citet[Theorem 3.1]{StSc12a} show that 
the injectivity of $\Ikn:H\to \Ltn$ is equivalent to either of 
the following statements:
\begin{enumerate}
	\item $(\sqrt{\lb_i} e_i)_{i\in I}$ is an ONB of $H$.
	\item For all $x,x'\in X$ we have 
		\begin{align}\label{kernel_sum_cont}
   k(x,x') = \sum_{i\in I} \lb_i e_i(x) e_i(x')\, .
\end{align}
\end{enumerate}
Obviously, if one of these conditions is true, then $H$ is separable, and 
 \citet[Corollary 2.10]{StSc12a} show that  the convergence in \eqref{kernel_sum_cont} is absolute
and we even  have $k(\cdot, x) = \sum_{i\in I}\lb_i e_i(x) e_i$ with unconditional convergence in $H$.


Let us now recall some notions related to  measures on topological spaces, see e.g.~\citet[Chapter IV]{Bauer01}
for details. 
To this end, let $(X,\t)$ be a Hausdorff (topological) space and $\nu$ be a  measure on its Borel-$\s$-algebra $\ca B(X)$.
Then $\nu$ is a Borel measure, if $\nu(K) < \infty$ for all compact $K\subset X$ and 
 $\nu$ is called strictly positive, if $\nu(O) > 0$ for all non-empty $O\in \t$.
Moreover, a finite  measure
$\nu$  on $\ca B(X)$ is a (finite) Radon measure if it 
is regular, i.e.~if for all $B\in \ca B(X)$ we have 
\begin{align*}
   \nu(B) & = \sup\{ \nu(K): K\mbox{ compact  and } K\subset B\} 
		= \inf\{ \nu(O): O\mbox{ open  and } B\subset O\}\, .
\end{align*}
A finite signed Radon measure is simply the difference of two finite Radon measures. 
In the following, we denote the space of all finite signed Radon measures by $\sosrm X$  
and the cone of (non-negative) finite Radon measures by $\sorm X$. As usual, 
$\sosrm X$ is equipped with the norm of total variation. 
%
%
Obviously, every finite Radon measure is a finite 
Borel measure, and by Ulam's theorem, see e.g.~\citet[Lemma 26.2]{Bauer01},
the converse implication is true if $X$ is a Polish space. In this respect recall that compact, metrizable spaces are 
Polish.

Now let $X$ be a locally compact Hausdorff space and $C_0(X)$ be the space of 
continuous functions vanishing at infinity. As usual, we equip 
$C_0(X)$ with the supremum norm.  Then, Riesz's representation theorem 
for locally compact spaces, see e.g.~\citet[Theorem 20.48 together with Definition 20.41, Theorem 12.40, Definition 12.39,
and a simple translation into the real-valued case using 
Theorem 12.36]{HeSt65} shows that 
\begin{align}\nonumber
   \sosrm X &\to C_0(X)'\\ \label{riesz-repres}
			\mu &\mapsto \Bigl(f\mapsto \langle f,\mu\rangle := \int_Xf\diff\mu\Bigr)
\end{align}
is an isometric isomorphism. In the compact case, in which 
$C_0(X)$ coincides with the space of 
continuous functions $C(X)$,
this can also be found in 
 e.g.~\citet[p.~265, Theorem IV.6.3]{DuSc58}.

Given a locally compact Hausdorff space $X$, a continuous kernel $k$ on $X$ with RKHS $H$
is called universal if $H\subset C_0(X)$ and $H$ is dense in $C_0(X)$ with respect to $\inorm\cdot$.
Note that for compact $X$ the inclusion  $H\subset C_0(X)$ is automatically satisfied.
Examples of universal kernels as well as various necessary and sufficient conditions
for universality can be found in e.g.~\cite{Steinwart01a,MiXuZh06a,SriperumbFukumizuETAL2011,ChWaZh16a}
and the references mentioned therein.

\section{New Characterizations}

In this section we first compare the norms $\hnorm\cdot$ and $\tvnorm\cdot$
and show that  in infinite dimensional they are never equivalent. 
By establishing some structural result for characteristic kernels, we then 
demonstrate that characteristic kernels cannot reliably distinguish between
distributions that are far away with respect to $\tvnorm\cdot$.
We further relate the eigenfunctions of the kernel to the metric $\g_k$
and with the help of this relation we 
investigate continuous kernels on (locally) compact spaces. Finally, 
we characterize on which compact spaces $X$ there do exist characteristic kernels.

\subsection{General results}\label{sec:general}

In this subsection we investigate the semi-norm $\hnorm \cdot$ on 
$\sosbm X$ for bounded kernels and general $X$.
We begin with a result that compares $\hnorm \cdot$ with $\tvnorm\cdot $.

\begin{theorem}\label{norm-equivalence}
   Let $(X,\ca A)$ be a measurable space and $H$ be the RKHS of a bounded and measurable kernel $k$ on $X$
such  that the kernel embedding $\P:\sosbm X\to H$ is injective.
Then, the following statements are equivalent:
\begin{enumerate}
   \item The space $\sosbm X$ is finite dimensional.
	\item The norms $\hnorm \cdot$ and $\tvnorm\cdot$ on $\sosbm X$ are equivalent.
	\item The norm $\hnorm \cdot$ on $\sosbm X$ is complete.
	\item The kernel embedding $\P:\sosbm X\to H$ is surjective.
\end{enumerate}
\end{theorem}

Theorem \ref{norm-equivalence} shows that for most cases of interest $(\sosbm X,\hnorm\cdot)$ is \emph{not} a Hilbert space.
To illustrate the fourth statement of Theorem \ref{norm-equivalence}, recall that 
the space 
\begin{displaymath}
   H_{\mathrm{pre}} :=  \biggl\{  \P(\mu):   \mu \in \spann\{\d_x: x\in X  \}   \biggr\}
\end{displaymath}
 is dense in $H$, see  e.g.~\citet[Theorem 4.21]{StCh08}.
Moreover, the space  $\spann\{\d_x: x\in X  \}$ is, in a weak sense, dense in $\ca M(X)$,
and therefore it is natural to ask whether every $f\in H$ is of the form $f=\P(\mu)$ for some $\mu \in \sosbm X$.
Theorem \ref{norm-equivalence} tells us that the answer is no, unless 
$\sosbm X$ is finite dimensional. 
In this respect recall that   it has been recently mentioned by  \citet{SGScXXa} that the 
kernel embedding $\P$ is, in general, not surjective. However, the authors do not provide
any example, or conditions, for non-surjective $\P$.
Two \emph{examples} 
of non-surjective kernel embeddings $\P$ are provided by \citet[Section 3]{PiWuLiMuWo07a}, while our
 Theorem \ref{norm-equivalence} shows that 
actually \emph{all} injective $\P$ fail to be surjective whenever we have $\dim \sosbm X = \infty$.

Our next goal is to show that for characteristic kernels on infinite dimensional spaces $\sosbm X$
there always exists probability measures that have maximal $\tvnorm\cdot$-distance  but
arbitrarily small $\hnorm\cdot$-distance. To this end, we need a couple of preparatory results.
We begin with the following lemma that investigates the effect of $\eins_X\in H$.

\begin{lemma}\label{m0-props}
  Let $(X,\ca A)$ be a measurable space and $H$ be the RKHS of a bounded and measurable kernel $k$ on $X$.
  If  $\eins_X\in H$, then $\sosbmh 0 {} X$ is $\hnorm\cdot$-closed in $\sosbm X$, and if $k$ is, in addition,
  characteristic, then the kernel embedding $\P:\sosbm X\to H$ is injective.
\end{lemma}

The next simple lemma computes the $\hnorm\cdot$-norm of measures if $H$ is an RKHS 
of the sum of two kernels.

\begin{lemma}\label{char-add}
  Let $(X,\ca A)$ be a measurable space, and $k_1$, $k_2$ be 
    bounded measurable kernels on $X$ with RKHSs $H_1$ and $H_2$. Let $H$ be the RKHS of the kernel
    $k=k_1+k_2$. Then for all $\mu\in \sosbm X$ we have 
    \begin{displaymath}
     \hnorm \mu^2 = \snorm \mu_{H_1}^2 + \snorm \mu_{H_2}^2\, .
    \end{displaymath}
    In particular, if $k_1$ is characteristic or has an injective kernel embedding, then the same is true for $k$.
\end{lemma}

In \citet[Corollary 11]{SriperumbGrettonETAL2010} it has already be show that 
the sum of two bounded, continuous translation-invariant kernels on $\R^d$ is characteristic, if 
at least one summand is characteristic. Lemma \ref{char-add} shows that this kind of 
inheritance holds in the general case.

Our next lemma considers products of kernels. In particular it shows that such 
products can only be characteristic if the involved factors are strictly integrally positive definite.

\begin{lemma}\label{product-form}
    Let $(X_1,\ca A_1)$ and $(X_2,\ca A_2)$ be  
    measurable spaces and $k_1$, $k_2$ be a bounded, measurable kernels on $X_1$ and $X_2$,
    respectively. We denote the RKHSs of $k_1$ and $k_2$ by $H_1$ and $H_2$. Moreover, let $H$ be the 
	RKHS of the kernel $k:= k_1 \cdot k_2$ on $X_1\times X_2$. Then, for all $\mu_1\in \sosbm {X_1}$
	and $\mu_2\in \sosbm{X_2}$ we have 
	\begin{displaymath}
	   \hnorm{\mu_1\otimes \mu_2} = \snorm{\mu_1}_{H_1} \cdot \snorm{\mu_2}_{H_2}\, .
	\end{displaymath}
	In particular, if $\dim \sosbm{X_1} \geq 2$ and $\dim \sosbm{X_2}\geq 2$, and $k$ is characteristic, then 
	$k_1$ and $k_2$ are strictly integrally positive definite with respect to $\sosbm {X_1}$ and $\sosbm {X_2}$,
	respectively.
\end{lemma}

At first glance it seems that Lemma \ref{product-form} contradicts \citet[Corollary 11]{SriperumbGrettonETAL2010},
which shows that the product of two bounded, continuous, translation-invariant kernels on $\Rd$ is 
characteristic on $\Rd$ as soon as at least one factor is characteristic. However, a closer look reveals 
that their result
considers the restriction of the product to the diagonal, whereas we treat the unrestricted kernel.
Later in Corollary \ref{product-char}, we will see that, on compact spaces, the product of two 
strictly integrally positive definite kernels gives a strictly integrally positive definite kernel on the product.

The following lemma compares strictly integrally positive definite kernels with respect
$\sosbm X$ and $\sosbmh 0 {}X$. In an implicit form, it has already been  used, 
and a similar statement is \citet[Theorem 32]{SGScXXa}.

\begin{lemma}\label{plusone-sipd}
   Let $(X,\ca A)$ be a measurable space, and $k$  be a
    bounded measurable kernel on $X$. Moreover, let $\ca M\subset \sosbm X$ be a subspace 
    with $\ca M\cap \sosbmh 1 {}X\neq \emptyset$
    and 
    $\ca M_0 := \ca M \cap \sosbmh 0 {} X$.
    Then the following statements are equivalent:
    \begin{enumerate}
     \item $k$ is strictly integrally positive definite with respect to $\ca M_0$.
     \item $k+1$ is strictly integrally positive definite with respect to $\ca M$.
     \item $k+1$ is strictly integrally positive definite with respect to $\ca M_0$.
    \end{enumerate}
\end{lemma}

We have already seen in Theorem \ref{norm-equivalence} 
that in the infinite-dimensional case the  norms $\hnorm\cdot$ and $\tvnorm\cdot$ are not 
equivalent on $\sosbm X$. Intuitively, this carries over to the subspace $\sosbmh 0{}X$.
The following result confirms this intuition as long as $\sosbmh 0{}X$ is a $\hnorm\cdot$-closed
subspace of $\sosbm X$.

\begin{theorem}\label{non-equivalence}
   Let $(X,\ca A)$ be a measurable space such that 
  $\dim \sosbm X = \infty$ and 
  $H$ be the RKHS of a bounded and measurable kernel $k$ on $X$
such  that the kernel embedding $\P:\sosbm X\to H$ is injective.
If $\sosbmh 0{}X$ is a $\hnorm\cdot$-closed
subspace of $\sosbm X$,
then  $\hnorm\cdot$ and $\tvnorm\cdot$ are not equivalent on 
  $\sosbmh 0 {} X$. 
\end{theorem}

The non-equivalence of $\hnorm \cdot$ and $\tvnorm\cdot$ on $\sosbmh 0{} X$ has already been observed in 
some particular situations. For example, 
\citet[Theorem 23]{SriperumbGrettonETAL2010} show that for universal kernels on compact metric spaces, 
$\g_k$ metrizes the weak topology (in probabilist's terminology) 
on $\sobpm X$, and since for $\dim \sosbm X = \infty$
this weak topology is strictly coarser than the $\tvnorm\cdot$-topology, we see that 
$\hnorm \cdot$ and $\tvnorm\cdot$ cannot be equivalent for such kernels. 
In addition, the non-equivalence can also 
be obtained from  
\citet[Theorems 21 and 24]{SriperumbGrettonETAL2010} for other continuous kernels on certain  metric spaces.
Finally recall that 
$\sosbmh 0{}X$ is a $\hnorm\cdot$-closed
subspace of $\sosbm X$ if $\eins_X\in H$ by Lemma \ref{m0-props}.

With the help of Theorem \ref{non-equivalence} the next result shows that 
characteristic kernels cannot reliably distinguish between distributions that are far away in total variation norm.

\begin{theorem}\label{different-distances}
  Let $(X,\ca A)$ be a measurable space such that 
  $\dim \sosbm X = \infty$  
  and $H$ be the RKHS of a characteristic kernel $k$ on $X$.
  Then 
  for all $\e>0$ there exist distributions $Q_1,Q_2\in \sosbmh 1 {} X$ such that 
  $\tvnorm{Q_1-Q_2} = 2$ and $\hnorm{Q_1-Q_2}\leq \e$.
\end{theorem}

Theorem \ref{different-distances} only shows that there are some distributions that cannot 
be reliably distinguished. The following corollary 
shows that such distributions actually occur everywhere.

\begin{corollary}\label{different-distances-cor}
  Let $(X,\ca A)$ be a measurable space such that 
  $\dim \sosbm X = \infty$  
  and $H$ be the RKHS of a characteristic kernel $k$ on $X$.
  Then for all $P\in \sosbmh 1 {} X$, $\d\in (0,2]$, and  $\e\in (0,\d)$ there 
	  exist   $Q_1,Q_2\in \sosbmh 1 {} X$ such that $\tvnorm{P-Q_i}\leq \d$ for $i=1,2$, 
  $\tvnorm{Q_1-Q_2} = \d$, and $\hnorm{Q_1-Q_2}\leq \e$.
\end{corollary}


The next goal of this subsection is to investigate the $\hnorm\cdot$-norm 
with the help of the eigenvalues and -functions of 
the integral operator $T_{k,\nu}$.
We begin with the following lemma that computes the  inner product \eqref{hinner}
by these  eigenvalues and -functions.


\begin{lemma}\label{square-integral-formula}
Let $(X,\ca A, \nu)$ be $\s$-finite measure space and $k$ be a bounded, measurable kernel
with RKHS $H$
 for which 
$\Ikn:H\to \Ltn$ is compact and injective. Then, for all $\mu_1, \mu_2 \in \sosbm X$, we have 
\begin{displaymath}
   \int_X\int_X k(x,x') \, \diff\mu_1(x)\diff\mu_2(x') = \sum_{i\in I} \lb_i \cdot \Bigl(\int_X e_i \, \diff\mu_1  \Bigr) \cdot \Bigl(\int_X e_i \, \diff \mu_2  \Bigr)\, ,
\end{displaymath}
where $(\lb_i)_{i \in I}\subset (0,\infty)$ and $(e_i)_{i \in I} \subset H$ are as at \eqref{kernel_sum_cont}.
\end{lemma}

For an interpretation of this lemma, we write, for bounded measurable $f:X\to \R$ and $\mu\in \sosbm X$, 
\begin{displaymath}
   \langle f, \mu\rangle := \int_X f\, \diff\mu \, .
\end{displaymath}
Combining Lemma \ref{square-integral-formula} with \eqref{hinner}
we then have 
\begin{displaymath}
   \bigl\langle \P(\mu_1) ,  \P(\mu_2) \bigr\rangle_H = \sum_{i\in I} \lb_i  \langle e_i, \mu_1\rangle \langle e_i, \mu_2\rangle \, .
\end{displaymath}
In other words, all calculations regarding inner products and norms of the kernel embedding 
$\mu\mapsto \P(\mu)$ can be carried over to a weighted $\ell_2$-space.

To formulate the following theorem, we denote,
for a 
$\s$-finite measure $\nu$ on $(X,\ca A)$, the set of all $\nu$-probability densities contained 
in $\Ltn$ by $\D(\nu)$, that is
\begin{displaymath}
   \D(\nu):= \Bigl\{ \aec h\in \Ltn \cap \Lx 1 \nu: \aec h\geq 0 \mbox{ and } \int_X h\diff\nu = 1  \Bigr\}\, .
\end{displaymath}
Moreover, we write 
$\ca P_2(\nu) := \{ h \diff\nu: \aec h\in \D(\nu)\}$ for 
 the corresponding set of probability measures.

 With the help of these preparations we can now formulate the following 
 theorem that characterizes non-characteristic kernels on $\ca P_2(\nu)$
 and that also establishes a result similar to Theorem \ref{different-distances} 
 for non-characteristic kernels.

%
%

\begin{theorem}\label{kernel-metric-thm}
   Let $(X,\ca A, \nu)$ be $\s$-finite measure space and $k$ be a bounded, measurable kernel
with RKHS $H$
 for which 
$\Ikn:H\to \Ltn$ is compact and injective. 
Then for all $\aec h,\aec g \in \D(\nu)$ and $P:= h\diff\nu$, $Q:=g\diff\nu$ the kernel 
mean distance can be computed by 
\begin{equation}\label{kernel-metric-thm-hxx}
   \g_k^2(P,Q) = \sum_{i\in I} \lb_i  \langle \aec{h-g}, [e_i]_\sim\rangle_{\Ltn}^2\, .
\end{equation}
Moreover, the following statements are equivalent:
\begin{enumerate}
   \item There exist $Q_1,Q_2\in \ca P_2(\nu)$ with $Q_1\neq Q_2$ and $\g_k^2(Q_1,Q_2)=0$.
		\item There exists an $\aec f\in \Lx 1 \nu \cap [H]_\sim^\perp$ with $\aec f \neq 0$ and 
$\int_X f \diff\nu = 0$.
	\item There exist $\aec{h_1},\aec{h_2}\in \D(\nu)$ with $\aec {h_1}\neq \aec{h_2}$ such that for all $i\in I$ we have 
	\begin{displaymath}
	   \langle \aec{h_1}, [e_i]_\sim\rangle_{\Ltn} =  \langle \aec{h_2}, [e_i]_\sim\rangle_{\Ltn}\, .
	\end{displaymath}
\end{enumerate}
Moreover, if one, and thus all, statements are true we actually find 
%
%
for all $P\in \ca P_2(\nu)$ 
and $\e\in (0,2)$ 
some $Q_1,Q_2\in \ca P_2(\nu)$ with 
$\tvnorm {P-Q_i}\leq \e$, $\tvnorm{Q_1-Q_2}=\e$, and 
 $\g_k^2(Q_1,Q_2)=0$. 
\end{theorem}

Equation \eqref{kernel-metric-thm-hxx} can also be used to 
show that under certain circumstances $\hnorm\cdot$ cannot reliably identify, for example, 
the uniform distribution. The following result, which is particularly interesting in view of 
Section \ref{sec:structure}, illustrates this.

\begin{corollary}\label{no-uniform}
      Let $(X,\ca A, \nu)$ be a probability space and $k$ be a bounded, measurable kernel
with RKHS $H$
 for which 
$\Ikn:H\to \Ltn$ is compact and injective. Assume that there is one eigenfunction $e_{i_0}$ with 
 $e_{i_0} = \eins_X$. In addition assume that  there are constants $c_1>0$ and $c_\infty<\infty$ with 
$\snorm{e_i}_{\Lx 1 \nu}\geq c_1$ and $\inorm{e_i}\leq c_\infty$ for all $i\in I$.
For $\a := c_\infty^{-1}$ and $j\neq i_0$ consider the signed measure
$Q_{j} := (\eins_X + \a e_j) d\nu$. Then $Q_{j}$ is actually a probability measure and for $P:= \nu$ we have 
\begin{align*}
   \tvnorm{P-Q_j} &\geq  c_1  c_\infty^{-1} \\
	\hnorm{P-Q_j} & = \lb_j  c_\infty^{-2}\, .
\end{align*}
\end{corollary}

The last result of this subsection provides some necessary conditions for characteristic kernels.

\begin{corollary}\label{cor:codim}
   Let $(X,\ca A, \nu)$ be finite measure space and $k$ be a bounded, measurable kernel
with RKHS $H$
 for which 
$\Ikn:H\to \Ltn$ is compact and injective. 
Then the following statements are true:
\begin{enumerate}
\item If $\codim [H]_\sim \geq 2$ in $\Ltn$, then $k$ is not characteristic.  
 \item If $\codim [H]_\sim \geq 1$ in $\Ltn$ and $\eins_X\in H$, then $k$ is not characteristic.
\end{enumerate}
\end{corollary}

\subsection{Continuous Kernels on Locally Compact Subsets}\label{sec:lcs}

In this subsection, we apply the general theory developed so far to 
bounded continuous kernels on locally compact Hausdorff spaces $(X,\t)$.

Let us begin with some preparatory remarks. To this end, let 
 $k$ be a bounded and continuous kernel on $X$ whose RKHS $H$ satisfies 
 $H\subset C_0(X)$. In the following we call such a $k$ a $C_0(X)$-kernel.
Our goal in this section is to investigate when $C_0(X)$-kernels are universal or 
characteristic. We begin with the following result that provides a necessary 
condition for the existence of strictly integrally positive definite kernels.

\begin{theorem}\label{exist-sipd}
 Let $(X,\t)$ be a locally compact Hausdorff space with $\sosrm X \neq \sosbm X$. Then 
 no $C_0(X)$-kernel  is strictly integrally positive definite with respect to $\sosbm X$.
\end{theorem}

Note that \citet{SriperumbFukumizuETAL2011} restrict their considerations to characteristic kernels
on locally compact \emph{Polish} spaces, for which we automatically have $\sosrm X = \sosbm X$
by Ulam's theorem, see e.g.~\citet[Lemma 26.2]{Bauer01}.
Some other papers, however, do not carefully distinguish between 
Borel and Radon measures, which, at last consequence, means that their results only hold if
we additionally assume $\sosrm X = \sosbm X$. 
Theorem \ref{exist-sipd} shows that this restriction is natural, and actually no restriction at all.
Furthermore note that 
for compact spaces $X$ one can use  Theorem \ref{exist-sipd} to show that 
$\sosrm X = \sosbm X$ is necessary for the existence of characteristic kernels.
We skip such a result since later
in Theorem \ref{exists-univer-char-kern},
we will be able to show an even stronger result.

%
%

Before we formulate our next result we need a bit more preparation. To this end, let 
$k$ be a $C_0(X)$-kernel on a locally compact space $(X,\t)$.
 Then we have $H\subset C_0(x)$ and 
a quick closed-graph argument shows 
that the corresponding inclusion 
operator $I:H\to C_0(X)$ is bounded. By the identification $C_0(X)'= \sosrm X$ in 
\eqref{riesz-repres} and the simple calculation 
\begin{displaymath}
   \langle Ih, \mu\rangle_{C_0(X),\sosrm X} = \int_X h\diff\mu  = \int_X \langle h, k(x,\cdot)\rangle_H\diff\mu(x)
= \langle h, \P(\mu)\rangle_H\, ,
\end{displaymath}
which holds for all $h\in H$, $\mu\in \sosrm X$ we further find that the adjoint $I'$ of $I$ is 
given by $I'= \P$.
This simple observation leads to 
the following characterization, which has already been shown for 
compact spaces $X$ by \citet[Proposition 1]{MiXuZh06a} and for locally compact Polish spaces 
by  \citet[Proposition 4]{SriperumbFukumizuETAL2011}. Although the proof of the latter paper also works 
on general locally compact Hausdorff spaces,  we decided to 
add the few lines for the sake of completeness.

\begin{theorem}\label{thm-new-char-2} 
   Let $(X,\t)$ be a locally compact Hausdorff space
and  $k$ be a $C_0(X)$-kernel.
Then the following two  statements are equivalent:  
\begin{enumerate}
   \item $k$ is strictly integrally positive definite with respect to $\sosrm X$.
	\item $k$ is universal.
\end{enumerate}
\end{theorem}


With the help of Theorem \ref{thm-new-char-2} we can now show that for characteristic kernels on
compact spaces $X$
it suffices to consider 
metrizable $X$. A similar result for universal kernels, which is included in the following theorem, has already been 
derived by \citet{StHuSc06a}.

\begin{theorem}\label{exists-univer-char-kern}
For a compact topological Hausdorff space $(X,\t)$ the following statements are equivalent:
\begin{enumerate}
\item There exists a universal kernel $k$ on $X$.
\item There exists a continuous characteristic kernel $k$ on $X$.
\item $X$ is metrizable, i.e.~there exists a metric generating the topology $\t$.
\end{enumerate}
If one and thus all statements are true, $(X,\t)$ is a compact Polish space and $\sosrm X = \sosbm X$.
\end{theorem}

Theorem \ref{exists-univer-char-kern} shows that on compact spaces 
 we may only expect universal or characteristic
kernels, if the topology is metrizable. Since in this case we have $\sosrm X = \sosbm X$, 
Theorem \ref{thm-new-char-2} and Proposition \ref{M0-char}
show the well-known result that every universal kernel  is characteristic.
In general, the converse implication is not true, but adding some structural requirements,
both notions may coincide. 
The following corollary illustrates this by showing that for product kernels universal and characteristic kernels 
coincide. 
 
\begin{corollary}\label{product-char}
    Let $(X_1,\t_1)$ and $(X_2,\t_2)$ be non-trivial compact metrizable spaces 
      and $k_1$, $k_2$ be continuous  kernels on $X_1$ and $X_2$,
    respectively. For the kernel $k:= k_1 \cdot k_2$ on $X_1\times X_2$ the following statements 
	are then equivalent:
	\begin{enumerate}
	   \item $k$ is universal.
		\item $k$ is characteristic.
		\item $k_1$ and  $k_2$ are universal.
	\end{enumerate}
\end{corollary}

Our next theorem, which provides a characterization of universal kernels 
with the help of the eigenfunctions of
the integral operator 
 $T_{k,\nu}$, is an extension of \citet[Corollary 5]{MiXuZh06a}
from compact to arbitrary locally compact Hausdorff spaces.
Before we present it, let us first make some preparatory remarks.
 To this end, let 
$\nu$ be a 
strictly positive and $\s$-finite Borel measure on $X$. For the matter of 
concreteness note that if  $X$ contains a dense, countable subset $(x_i)_{i\geq 1}$ then 
$\nu:=\sum_{i\geq 1} \d_{x_{i}}$ satisfies these assumptions and therefore we always find such measures on 
e.g.~compact metric spaces.
Now, let $k$ be a bounded and continuous kernel on $X$ satisfying 
\eqref{int-diag}. 
Then $H$ consists of continuous functions and 
 \citet[Corollary 3.5]{StSc12a}  show that \eqref{kernel_sum_cont} holds for all $x,x'\in X$,
and consequently, the assumptions of Lemma \ref{square-integral-formula},
Theorem \ref{kernel-metric-thm}, and  Corollary \ref{cor:codim} are satisfied.

With these preparations we can now formulate the following characterization of universal kernels,
where we note that the equivalence between \emph{i)} and \emph{ii)}
has essentially been shown in \cite[Proposition 12]{SrFuLa10a}.

\begin{theorem}\label{thm-new-char-3} 
   Let $(X,\t)$ be a locally compact Hausdorff space, $\nu$ be a strictly positive, $\s$-finite  Borel measure on $X$,
and  $k$ be a $C_0(X)$-kernel  satisfying \eqref{int-diag}.
In addition, let $(e_i)_{i\in I}$   be the eigenfunctions of $T_{k,\nu}$ in 
\eqref{kernel_sum_cont}.
Then the following statements are equivalent: 
\begin{enumerate}
	\item $k$ is universal.
	\item For all $\mu\in \sosrm X$ satisfying $\int_X e_i \diff\mu =0$ for all $i\in I$ 
	we have $\mu=0$.
	\item The space $\spann\{ e_i: i\in I\}$ is dense in $C_0(X)$.
\end{enumerate}
If one, and thus all, statements are true and $\nu\in \sosbmh {}*X$, then
 $([e_i]_\sim)_{i\in I}$ is an ONB of $\Ltn$.
\end{theorem}

Our next result characterizes universal and characteristic kernels on compact spaces 
with the help of the eigenfunctions and -values of a suitable $T_{k,\nu}$.
In view of  Theorem \ref{exists-univer-char-kern} it suffices to consider  compact spaces that are
 Polish.

\begin{corollary}\label{uni-cor}
Let $(X,\t)$ be a compact  metrizable space 
and  $k$ be a continuous kernel   with RKHS $H$. 
 Moreover, let $(\lb_i)_{i\in I}\subset [0,\infty)$ be a family  converging to $0$ and 
$(e_i)_{i\in I} \subset C(X)$  be a family
such that $\spann\{e_i:i\in I\} $ is dense in $C(X)$ and
\begin{displaymath}
 k(x,x') = \sum_{i\in I} \lb_i e_i(x)e_i(x')
\end{displaymath}
holds for all $x,x'\in X$. If there is a 
strictly positive, finite and regular Borel measure $\nu$  on $X$ such that 
$([e_i]_\sim)_{i\in I}$
is an ONS in $\Ltn$, then:
\begin{enumerate}
   \item $k$ is universal if and only if $\lb_i>0$ for all $i\in I$. 
		\item If $e_{i_0}=\eins_X$ for some $i_0\in I$, then $k$ is characteristic if and only if 
	$\lb_i >0$ for all $i\neq i_0$.
		\item If $\eins_X \in H$ and $e_{i}\not=\eins_X$ for all $i\in I$, then $k$ is characteristic if and only if 
	$\lb_i >0$ for all $i\in I$.
\end{enumerate}
\end{corollary}

\section{Characteristic kernels on spaces with additional structure}\label{sec:structure}

In this section, we apply the developed theory to translation-invariant or isotropic kernels on compact Abelian groups or 
spheres, respectively.


\subsection{Compact Abelian Groups}\label{sec:groups}

In this subsection we apply the theory developed so far to translation-invariant 
kernels on compact Abelian groups. Here the main difficulty lies in the fact that 
one traditionally  considers  kernels on groups that are $\C$-valued, while we are only interested 
in $\R$-valued kernels. Although at first glance, one may not expect  any problem
arising from this discrepancy,
it turns out that it actually does make a difference when constructing an
ONB of $\Ltn$ with the help of characters as soon as we have more than one self-inverse 
character.

Our first goal is to make the introducing remarks precise. To this end, 
 let $(G, +)$ be a compact Abelian group, and $\nu$ be its normalized Haar measure.
As usual, we write $\Lx 2 G := \Ltn$, and 
 for later use recall that $\nu$ is strictly positive and regular, see e.g.~\citet[p.~193/4]{HeRo63}.
Moreover, let 
 $(\hat G, \cdot)$ be the  dual group of $G$, which  consists of characters $e: G\to \T$,
where $\T$ denotes, as usual, the unit circle in $\C$, see e.g.~\citet[Chapter Six]{HeRo63}
and \citet[Chapter 4.1]{Folland95}.
 For notational purposes, we assume that we have another group $(I,+)$ that is isomorphic
 to $(\hat G, \cdot)$ by some mapping  $i\mapsto e_i$. This gives $e_{i+j} = e_i e_j$, $e_0 = \eins_G$,
 and since we further have $e \bar e= \eins_G$ for all $e\in \hat G$, our  notation also yields
 $e_{-i} =  \bar e_i$. In particular, we have $\re e_{-i} = \re e_i$ and $\im e_{-i} = -\im e_{i}$
for all $i\in I$, and for all $i\in I$ with $i=-i$ the latter equality immediately yields $\im e_i = 0$.
Finally, for $i\in I$ and $x,y\in G$ we have 
\begin{displaymath}
   e_i(-y+x) = \frac {e_i(x)}{e_i(y)} = \overline{e_i(y)} e_i(x) = e_{-i}(y) e_i(x)\, ,
\end{displaymath}
and from this it is easy to derive both, $\re e_i(-x) = \re e_i(x)$ and $\im e_i(-x) = -\im e_i(x)$,
as well as
\begin{align}\label{add-thm}
\re e_i(-y+x) = \re e_i(x) \re e_i(y) + \im e_i(x)\im e_i(y)\, .
\end{align}
Note that for $i\in I$ with $i=-i$, the latter formula can be simplified using $\im e_i = 0$.

Now, it is well-know that 
 $([e_i]_\sim)_{i\in I}$
 is an ONB of $\Lx 2 {G,\C}$, see e.g.~\citet[Corollary 4.26]{Folland95}
 and using this fact a quick application of the Stone-Weierstrass theorem shows that 
 $(e_i)_{i\in I}$ 
 is also dense in $C(G, \C)$. Let us construct a corresponding ONB in $\Lx 2 G$. To this end,
we write $I_0 := \{i\in I: i=-i\}$
for the set of all \emph{self-inverse} elements of $I$.
 Moreover, we fix a partition $I_+\cup I_- = I\setminus I_0$ 
such that $i\in I_+$ implies $-i\in I_-$ for all $i \in I\setminus I_0$.
Obviously, the sets $I_0, I_+, I_-$ form a partition of $I$. 
Let us now define the family $(e_i^*)_{i\in I}$ by
\begin{equation}\label{real-onb-def}
   e_i^* :=  
\begin{cases}
 \re e_i & \mbox{ if }\, i\in I_0 \\
 \sqrt{2}  \re e_i & \mbox{ if }\, i\in I_+ \\
	\sqrt{2} \im e_i  & \mbox{ if }\, i\in I_- \, .
\end{cases}
\end{equation}
The next result shows that $(e_i^*)_{i\in I}$ is the desired family.

\begin{lemma}\label{real-onb}
 Let $(G, +)$ be a compact Abelian metric group. 
 Then each family
  $(e_i^*)_{i\in I}$ given by \eqref{real-onb-def}
  is an ONB of $\Lx 2 G$ and $\spann\{e_i^*: i\in I\}$ is dense in $C(G)$.
  Finally, we have $\inorm{e_i^*}\leq \sqrt 2$ for all $i\in I$.
\end{lemma}

In the following, we call a kernel $k$ on an Abelian group $(G,+)$
translation invariant, if there exists a 
function $\k:G\to \R$ such that $k(x,x') = \k(-x +x')$ for all $x,x'\in G$. 
Clearly, $k$ is continuous, if and only if $\k$ is. The following
lemma provides a representation of translation invariant kernels.

\begin{lemma}\label{real-bochner}
 Let $(G, +)$ be a compact Abelian group, 
  $(e_i^*)_{i\in I}$ be a family of the form \eqref{real-onb-def},
 and $k:G\times G\to \R$ be a bounded, measurable function. 
Then the following statements are equivalent:
\begin{enumerate}
   \item $k$ is a bounded, measurable, and translation invariant kernel on $G$.
	\item There exists a summable
 family $(\lb_i)_{i\in I}\subset [0,\infty)$ such that 
 \begin{equation}\label{Mercer-on-G}
     k(x,x')
 = \sum_{\lb_i > 0} \lb_i e_i^*(x) e_i^*(x')  = \sum_{\lb_i> 0} \lb_i \re e_i(x-x')\, ,
 \end{equation}
  where the series converge absolutely for all $x,x'\in G$ as well as uniformly in $(x,x')$. 
\end{enumerate}
If one, and thus both, statements are true, then $k$ is continuous, $(\lb_i)_{i\in I}$ are 
all, possibly vanishing, eigenvalues of $T_{k,\nu}$, and $(e_i^*)_{i\in I}$ is an ONB of the corresponding eigenfunctions.
%
\end{lemma}

For an interpretation of the representation \eqref{Mercer-on-G} 
 recall that $\k(-x+x') = k(x,x')$ for all $x,x\in G$ and hence  \eqref{Mercer-on-G} 
gives 
\begin{displaymath}
   \k(x) =  \sum_{\lb_i> 0} \lb_i \re e_i(-x) =  \sum_{\lb_i> 0} \lb_i \re e_i(x) 
\end{displaymath}
for all $x\in G$. Consequently, the second equality in \eqref{Mercer-on-G} 
is  Bochner's theorem, see e.g.~\citet[Theorem 30.3]{HeRo70}, in the case of 
$\R$-valued kernels on compact Abelian groups. Unlike this classical theorem, however, 
the second equality in \eqref{Mercer-on-G}  also describes 
 how the representing measure of $\k$ on $\hat G$ is given by the eigenvalues of 
$T_k$ or $T_k^\C$. In the following, we do not need this information, in fact, we only 
mentioned the second equality to provide a link to the existing theory.
Instead, the 
first equality in \eqref{Mercer-on-G}, which replaces the characters of $G$ by the eigenfunctions
of $T_k$, is more important for us, since this equality 
actually is the Mercer representation of $k$ in the sense of \eqref{kernel_sum_cont}
and therefore the theory developed in the previous sections becomes applicable.

The next result is an extension of Theorem \ref{exists-univer-char-kern} to translation-invariant 
kernels on compact Abelian groups.

\begin{theorem}\label{exist-univ-on-G}
   Let $(G,+)$ be a compact Abelian group. Then the following statements are equivalent:
\begin{enumerate}
   \item $G$ is metrizable.
	\item $\hat G$ is at most countable.
	\item There exists a translation-invariant universal kernel on $G$.
	\item There exists a universal kernel on $G$.
	\item There exists a translation-invariant  characteristic kernel on $G$.
	\item There exists a continuous characteristic kernel on $G$.
\end{enumerate}
\end{theorem}

Note that the equivalence between
\emph{i)} and \emph{ii)} can   also be shown without using translation-invariant kernels, 
see e.g.~\citet[Proposition 3]{Morris79a} or \citet[Theorem 24.15]{HeRo63}.
Our proof, to the contrary, is solely RKHS-based.

Our next result characterizes universal and characteristic translation-invariant kernels 
on compact Abelian groups. In view of Theorem \ref{exist-univ-on-G}, it suffices to consider 
the metrizable case.

\begin{corollary}\label{char-on-G-char}
 Let $(G,+)$ be a compact metrizable Abelian group
 and $k$ be a   translation-invariant
 kernel on $G$ with representation \eqref{Mercer-on-G}. Then we have:
 \begin{enumerate}
  \item $k$ is universal if and only if $\lb_i>0$ for all $i\in I$.
  \item $k$ is characteristic if and only if $\lb_i>0$ for all $i\neq 0$.
 \end{enumerate}
\end{corollary}

Corollary \ref{char-on-G-char} generalizes \citet[Theorem 14 and Corollary 15]{SriperumbGrettonETAL2010}
from $\T^d$ to arbitrary compact  metrizable Abelian groups. Moreover, recall that these authors also
provide a couple of translation-invariant characteristic kernels on $\T$ that 
enjoy a closed form.

%
%

As mentioned in the beginning of this section, the
major difficulty in deriving a Mercer representation \eqref{Mercer-on-G}
for translation-invariant kernels is the handling of self-inverse characters
other than the neutral element. The simplest example of a group $G$ whose dual $\hat G$ contains
more than one self-inverse is the quotient group $(\Z_2, +)$ of $(\Z,+)$ with its subgroup $2\Z$.
Indeed, besides the neutral element $e_0$, $\hat \Z_2$ only contains the 
character $e_1$ given by $e_1(0) := 1$ and $e_1(1) := -1$. 
Clearly, this gives $e_1^2 = e_0$ and thus $e_1$ is self-inverse.
Now note that a function
$k:\Z_2\times \Z_2\to \R$ can be uniquely described by a 2-by-2 matrix 
$K = (k(x,x'))_{x,x'\in \Z_2}$ and a simple calculation shows that $k$ is a   kernel
if and only if $k(0,1) = k(1,0)$, $k(0,0)\geq 0$, $k(1,1) \geq 0$, and $k(0,0)k(1,1) \geq k^2(0,1)$. 
Moreover, $k$ is translation-invariant as soon as  it is constant on its diagonal, and in this case
the previous conditions reduce to
\begin{equation}\label{discrete-trala}
   k(0,1) = k(1,0) \qquad \qquad \mbox{ and } \qquad \qquad  k(0,0) = k(1,1) \geq |k(0,1)|\, .
\end{equation}
Now, let $k$ be a translation-invariant kernel on $\Z_2$ and 
$\lb_0, \lb_1\geq 0$ be the coefficients in \eqref{Mercer-on-G}. Then
a simple calculation shows that
 the describing matrix $K$ is given by 
\begin{displaymath}
K=
  \begin{pmatrix}
    \lb_0+\lb_1 & \lb_0-\lb_1 \\
   \lb_0-\lb_1 & \lb_0+\lb_1
  \end{pmatrix}\, ,
\end{displaymath}
and therefore  it is not hard to see by  Corollary \ref{char-on-G-char} 
that $k$ is characteristic, if and only if 
$k(0,0) \neq k(0,1)$. Similarly, $k$ is universal, if and only if 
$k(0,0) \neq \pm k(0,1)$.

While this example seems to be rather trivial, it already has some important applications.
For example, assume that our input space $X$ is a product space for which some components 
belong to a compact metrizable Abelian group, while the remaining components are only allowed to 
attain the values $0$ and $1$. In other words, $X$ is of the form 
$X = G \times \Z_2^d$, where $G$ is a compact metrizable Abelian group and $d\geq 1$.
Now, an intuitive way to construct a (translation-invariant)
 characteristic kernel $k$ on $X$ is to 
take a product $k:= k_C\cdot k_D$, where 
$k_C$ and $k_D$
are kernels  on $G$ and $\Z_2^d$, respectively. By Corollary \ref{product-char} we then 
know that $k$ is characteristic (or universal) if and only if both 
$k_C$ and $k_D$ are universal. Clearly, if 
$k_D$ is itself a product of  kernels $k_1,\dots,k_d$ 
then  $k_D$ is almost automatically translation-invariant and universal. Indeed, 
if all $k_i$
satisfy
\eqref{discrete-trala} with  $k_i(0,0) \neq \pm k_i(0,1)$, then each $k_i$ is translation-invariant and universal, 
and thus so is $k_D$.
It seems fair to say that most ``natural'' choices of $k_i$ will satisfy these assumptions.
On the other hand, translation-invariant universal kernels $k_C$ on $G$ are completely characterized by Corollary
\ref{char-on-G-char}, and thus it is straightforward to characterize all translation-invariant
characteristic kernels $k$ on $G\times \Z_2^d$ of product type $k:= k_C\cdot k_D$.
However, their representation \eqref{Mercer-on-G} is a bit more cumbersome. Indeed, any element
$(e, \om)\in \hat G \times \hat \Z_2^d = \widehat{G \times \Z_2^d}$
with self-inverse $e\in \hat G$ and arbitrary
$\om \in \Z_2^d$ is self-inverse, where the 
equality in the sense of group isomorphisms can e.g.~be found in 
\citet[Proposition 4.6]{Folland95}.
Consequently, the set $I_0$, which intuitively may be viewed as a small set of unusual characters,
may actually be rather large. 

Note that the set  $\Z_2$  appears 
in data analysis settings whenever we have categorical variables with two possible values,
which quite frequently is indeed the case.
Now we have seen around \eqref{discrete-trala} 
that most natural choices of kernels on $\Z_2$ are actually translation-invariant,
and for these kernels the results of this subsection applies. Similarly, if we have 
a categorical variable with an even number $m$ of possible values that have a cyclic nature,
for example hours of a day or months of a year, 
then $\Z_m$
can be used to describe these values, and kernels that respect the cyclic nature are
translation-invariant. Clearly $m/2$ is self-inverse in $\Z_m$, and therefore 
our theory again applies.

For more on structural properties of compact Abelian groups as well as for further examples 
 we refer to
\citet[\S 25]{HeRo63} and
 \citet[Chapter 8]{HoMo13}.

%

\subsection{Spheres}\label{sec:Sd}

In this subsection
we consider isotropic kernels $k$ on $\bbS^d$ for $d \ge 1$, that is
kernels of the form
 $k(x,y) = \psi(\theta(x,y))$, where $\theta(\cdot,\cdot)=\arccos \langle \cdot,\cdot\rangle$ is the geodesic distance on $\mathbb{S}^d$ and $\langle \cdot,\cdot\rangle$ denotes the standard scalar product on $\mathbb{R}^{d+1}$.   
Following \citet{Gneiting2013}, let $\Psi_d$ be the class of all continuous 
functions $\psi$ on $[0,\pi]$ such that $k(x,y) = \psi(\theta(x,y))$ is a kernel on $\bbS^d$, 
and define $\Psi_{\infty} := \cap_d \Psi_d$. 
We write $\Psi_d^+\subset \Psi_d$ 
for the class of functions that induce a strictly positive definite 
kernel on $\bbS^d$, and set $\Psi_{\infty}^+:= \cap_d \Psi_d^+$. It holds that $\Psi_{d+1} \subset \Psi_d$ and $\Psi_{d+1}^+ \subset \Psi_d^+$, see \citet[Corollary 1]{Gneiting2013}.

The following two theorems are our main results on characteristic kernels on $\bbS^d$. 

\begin{theorem}\label{thm:32} Let $k$ be an isotropic kernel on $\mathbb{S}^d$ induced by
some  $\psi \in \Psi_{d+2}$ or by some $\psi\in\Psi_{d+1}^+$. Then the following statements are equivalent: 
\begin{enumerate}
\item $k$ is characteristic. 
\item $k$ is universal.
\item $k$ is strictly positive definite on $\bbS^d$, that is, $\psi \in \Psi_d^+$.
\end{enumerate}
\end{theorem}

Theorem \ref{thm:32} shows in particular that any $\psi \in \Psi_{d+1}^+$ induces a characteristic kernel on $\mathbb{S}^d$. For the practically most relevant case of $\mathbb{S}^2$, all the parametric families of isotropic positive definite functions in \citet[Table 1]{Gneiting2013} are in $\Psi_3^+$ and thus all of them are characteristic and yield strictly proper kernel scores on $\mathbb{S}^2$ by Theorem \ref{thm:one} and Theorem \ref{thm:32}.

\begin{theorem}\label{thm:33}
Let $k$ be an isotropic kernel on $\mathbb{S}^d$ induced by $\psi \in \Psi_{\infty}$. Then the following statements are equivalent: 
\begin{enumerate}
\item $k$ is characteristic. 
\item $k$ is universal.
\item $k$ is strictly positive definite, that is, $\psi \in \Psi_{\infty}^+$.
\end{enumerate}
\end{theorem}


Theorem \ref{thm:33} is analogous to the result of \citet[Proposition 5]{SriperumbFukumizuETAL2011} for radial kernels on $\R^d$. 


For the proofs of Theorems \ref{thm:32} and \ref{thm:33}, we need to introduce some preliminaries.
By \citet{Schoenberg1942} the functions in $\Psi_{\infty}$ have a representation of the form
\begin{equation}\label{eq:Schoenberg}
\psi(\theta) = \sum_{n=0}^{\infty}b_n (\cos(\theta))^n\, , \qquad \theta \in [0,\pi],
\end{equation}
where $(b_n)_n$ is a summable sequence of non-negative coefficients termed the \emph{$\infty$-Schoenberg sequence} of $\psi$. 
\citet{Schoenberg1942} also showed that the functions in $\Psi_d$ have a representation as
\begin{equation*}
\psi(\theta) = \sum_{n=0}^{\infty}b_{n,d} \frac{C_n^{(d-1)/2}(\cos\theta)}{C_n^{(d-1)/2}(1)}\, , \qquad \theta \in [0,\pi],
\end{equation*}
where $(b_{n,d})_n$ is a summable sequence of non-negative coefficients termed the \emph{$d$-Schoenberg sequence} of $\psi$, and $C_n^\lambda$, $\lambda > 0$, $n \in \mathbb{N}_0$ are the Gegenbauer polynomials; see the \citet[18.3.1]{dlmf}. For $\lambda =0$, we set $C_n^0(\cos \theta) := \cos (n\theta)$. 

\begin{definition}
A sequence of non-negative real numbers $(b_n)_{n\in \mathbb{N}_0}$ fulfills \emph{condition $b$}, if $b_n > 0$ for  $\infty$-many even and $\infty$-many odd integers. 
\end{definition}

\begin{remark}\label{rem:spd2cb}
For $\psi \in \Psi_\infty$ or $\psi \in \Psi_d$, $d\ge 2$, the induced isotropic kernel is strictly positive definite if and only if its Schoenberg sequence fulfills condition $b$, see \citet{Menegatto1992,Menegatto1994} and \citet{ChenMenegattoETAL2003}. For $d = 1$, condition $b$ remains a necessary condition 
for $\psi \in \Psi_1^+$  as shown by  \citet{Menegatto1995}. However, it is not sufficient any more. A simple sufficient condition for $\psi \in \Psi_1^+$ that is useful for our purposes but which is not necessary is that there is $n_0$ such that $b_{n,1} > 0$ for all $n \ge n_0$. See \citet{MenegattoOliveiraETAL2006} for a necessary and sufficient condition in the case $d = 1$.
\end{remark}

\begin{lemma}\label{prop:1} 
If $\psi \in \Psi_{d+2}$, then it is strictly positive definite on $\bbS^d$ if and only if $b_{n,d}>0$ for all $n \in \mathbb{N}_0$.
\end{lemma}

The set $\Psi_d^+ \backslash \Psi_{d+2} \subset \Psi_d$ is not empty and also contains elements with $b_{n,d} > 0$ for all $n \ge 0$. To construct an explicit example, take any summable sequence $(b_n)_{n \in \mathbb{N}_0}$ of positive real numbers such that $b_{2} > (d(d+3)/2)b_0$. Let $\psi$ be the function with $d$-Schoenberg sequence $(b_n)_{n \in \mathbb{N}_0}$. Then $\psi \in \Psi_d^+$, fulfills $b_{n,d} > 0$ for all $n \ge 0$, and, by \citet[Corollary 4]{Gneiting2013}, it is not a member of $\Psi_{d+2}$.

\begin{lemma}\label{prop:34}
If $\psi \in \Psi_{d+1}^+$, then $b_{n,d} > 0$ for all $n \ge 0$.
\end{lemma}

After these preliminary considerations concerning Schoenberg sequences, we will now show Theorem \ref{thm:32} by applying Corollary \ref{uni-cor}. To this end, let $(e_{n,j})_{n \in \mathbb{N}_0,j=1,\dots,N(d,n)}$ 
denote an orthonormal basis of
spherical harmonics on $\bbS^d$. The polynomial $e_{n,j}$ has order $n$ and $N(d,n) = \binom{n+d}{n} - \binom{n+d-2}{n-2}$; see for example \citet[Theorem 3.1.4]{Groemer1996}, where we note that he works on $\bbS^{d-1}$ while we work on $\bbS^d$. In particular, $e_{0,0} = 1$. By \citet[Theorem 3.3.3]{Groemer1996},
\begin{equation}\label{eq:legendre}
\frac{C_n^{(d-1)/2}(\langle x,y\rangle)}{C_n^{(d-1)/2}(1)} = \frac{1}{N(d,n)}\sum_{j=1}^{N(d,n)}  e_{n,j}(x)e_{n,j}(y),
\end{equation}
hence any isotropic kernel on $\mathbb{S}^d$ induced by $\psi \in \Psi_d$ has a Mercer 
representation of the form  \eqref{kernel_sum_cont} with $\lambda_{n,j} = b_{n,d}/N(d,n)$. Moreover, \citet[Corollary 3.2.7]{Groemer1996} shows that the space $\spann\{e_{n,j}: n \in \mathbb{N}_0,j=1,\dots,N(d,n)\}$ is dense in $C(\mathbb{S}^d)$. 

Similar to \citet[Theorem 9]{SriperumbGrettonETAL2010} for translation invariant kernels on $\R^d$, Corollary \ref{uni-cor} yields the following theorem. 

\begin{theorem}\label{thm:Schoenberg}
The kernel induced by $\psi \in \Psi_d$ is 
\begin{enumerate}
\item universal if and only if $b_{n,d} > 0$ for all $n \ge 0$.
\item characteristic if and only if $b_{n,d} > 0$ for all $n \ge 1$. 
\end{enumerate}
\end{theorem}

For the proof of the converse of Theorem \ref{thm:33}, we need the following proposition. 

\begin{proposition}\label{prop:421}
Let $k$ be a kernel on $\mathbb{S}^d$  induced by $\psi \in \Psi_\infty$. If $k$ is characteristic, then $\psi$ is strictly positive definite, that is $\psi \in \Psi_\infty^+$. 
\end{proposition}

\section{Proofs}\label{sec:proofs}

\subsection{Proofs related to Section \ref{sec:intro}}

For the proof of Theorem \ref{thm:one} we assume that the general 
results on kernel mean embeddings recalled in Section \ref{sec:prelim}
up to Definition \ref{def-char-kern} of characteristic kernels are available.
Note that these results do not involve kernel scores at all, so that 
there is no danger of circular reasoning.

\vspace*{1ex}
\begin{proofof}{Theorem \ref{thm:one}} 
As explained around \eqref{square-finite}, the condition $P \in \mathcal{M}_1^k(X)$
 ensures that $\Phi(P)$  defined by \eqref{eq:inj} is indeed an element of $H$. 
Now \eqref{main-obs} follows from \eqref{hinner}, namely
%
\begin{align*}
\snorm{\P(P)-\P(Q)}_H^2 
&= \int\!\!\!\!\int\! k(x,y)\diff P(x)\diff P(y) + \int\!\!\!\!\int\! k(x,y)\diff Q(x)\diff Q(y) 
- 2\int\!\!\!\!\int\! k(x,y)\diff P(x)\diff Q(y)\\
&= 2 \left(\int S_k(P,x)\diff Q(x) - \int S_k(Q,x)\diff Q(x)\right).
\end{align*}
The remaining assertions directly follow from \eqref{main-obs}.
\end{proofof}

\subsection{Proofs related to Section \ref{sec:general}}

\begin{proofof}{Theorem \ref{norm-equivalence}}
   \atob {iv} {iii} By assumption, $\P:\sosbm X\to H$ is bijective and since $H$ is complete, so is 
$\hnorm \cdot$ on $\sosbm X$.

\atob {iii} {ii} Consider the identity map $\id: (\sosbm X, \tvnorm\cdot) \to (\sosbm X, \hnorm\cdot)$.
Since we have already seen in Section \ref{sec:prelim} that
\begin{displaymath}
\hnorm \mu 
\leq \mnorm{[\, \cdot\, ]_\sim :H\to \Lx 1 \mu}
\leq \inorm {k} \cdot \tvnorm \mu
\end{displaymath}
holds for all $\mu\in \sosbm X$, the identity map is continuous. In addition, it is, of course, bijective,
and since
both $(\sosbm X, \tvnorm\cdot)$ and  $(\sosbm X, \hnorm\cdot)$ are complete, the open mapping theorem,
see e.g.~\citet[Corollary 1.6.8]{Megginson98} shows that both norms are equivalent.

\atob {ii} i By the arguments of \citet[page 63f]{DiJaTo95}
the space $(\sosbm X, \tvnorm\cdot)$  is a so-called ${\cal L}_{1,\lb}$-space for all $\lb>1$, see 
\citet[page 60]{DiJaTo95} for a definition, while the Euclidean structure of $\hnorm\cdot$ on $\sosbm X$, which is
inherited from $(H,\hnorm\cdot)$,
shows that 
$(\sosbm X, \hnorm\cdot)$ is an ${\cal L}_{2,\lb}$-space for all $\lb>1$.
Let us now assume that 
$\dim  \sosbm X = \infty$.
Then, \citet[Corollary 11.7]{DiJaTo95}
shows that $(\sosbm X, \tvnorm\cdot)$ has only trivial type 1,
while $(\sosbm X, \hnorm\cdot)$ has optimal type 2.
However, the definition of the type of a space, see \citet[page 217]{DiJaTo95}, immediately shows that
equivalent norms always share their type, and hence 
 $\hnorm \cdot$ and $\tvnorm\cdot$ are not equivalent on $\sosbm X$.

\atob i {iv} If $\sosbm X$ is finite dimensional, the $\s$-algebra $\ca A$ is finite. Since $k$ is assumed to be measurable,
we then see that $k$ can only attain finitely many different values, and consequently, 
the canonical feature map $\P:X\to H$ also attains only finitely many different values, say 
$f_1,\dots, f_m\in H$. 
Using \citet[Theorem 4.21]{StCh08}, we conclude that $H = \spann\{f_1,\dots,f_m\}$.
We now define $A_i := \P^{-1}(\{f_i\})$ for $i=1,\dots,m$. By construction $A_1,\dots,A_m$ 
form a partition of $X$ with $A_i\in \ca A$ and $A_i\neq \emptyset$ for all  
$i=1,\dots,m$. Let us fix some $x_i\in A_i$.
For the corresponding Dirac measures we then have $\P(\d_{x_i}) = \P(x_i) = f_i$
and since $\P: \sosbm X\to H$ is linear we find $\spann\{f_1,\dots,f_m\} \subset \P(\sosbm X)$.
This shows the surjectivity of the kernel embedding.
\end{proofof}

\begin{proofof}{Lemma \ref{m0-props}}
 Let $(\mu_n) \subset \sosbmh 0 {} X$ be a sequence that converges to some $\mu\in \sosbm X$ in $\hnorm\cdot$.
 Then we have $\mu_n(X) = 0$ for all $n\geq 1$ and we need to show $\mu(X)=0$.
 The latter, however, follows from  \eqref{alternative-comp}, namely
 \begin{displaymath}
 |\mu(X)| =  \biggl| \int \eins_X \diff (\mu_n-\mu)  \biggr| \leq \hnorm{\eins_X} \hnorm{\mu_n-\mu} \to 0\, .
 \end{displaymath}
  Let us now assume that $k$ is characteristic. Then Proposition \ref{M0-char} 
%
shows that   $k$ is
strictly
integrally positive definite with respect to $\sosbmh 0 {} X$, and hence
$\hnorm{\mu} > 0$ for all $\mu\in \sosbmh 0 {} X\setminus\{0\}$.
Now let $\mu\in \sosbm X\setminus \sosbmh 0 {} X$ and $P$ be 
some probability measure on $X$. By \eqref{m0+p} we then find an $\a\in \R$ and some $\mu_0\in \sosbmh 0 {} X$
such that $\a P + \mu_0 = \mu$, and since $\mu\not\in \sosbmh 0 {} X$ we actually have $\a\neq 0$.
Using this decomposition, and  \eqref{alternative-comp} we now find
\begin{displaymath}
 \hnorm{\mu} = \hnorm{\a P + \mu_0} \geq \frac 1 {\hnorm{1_X}}  \biggr|\int \eins_X \diff (\a P + \mu_0)\biggr|
 = \frac {|\a|} {\hnorm{1_X}} > 0
\end{displaymath}
and hence  $\P:\sosbm X\to H$ is injective, see \eqref{sipd}.
\end{proofof}

\begin{proofof}{Lemma \ref{char-add}}
 By the definition of the $\hnorm\cdot$-norm and \eqref{hinner} we have 
 \begin{displaymath}
  \hnorm \mu^2 = \int_X \int_X k_1(x,x') + k_2(x,x') \diff\mu(x)\diff \mu(x')
  =  \snorm \mu_{H_1}^2 + \snorm \mu_{H_2}^2\, .
 \end{displaymath}
 Moreover, if $k_1$ is characteristic, then $\snorm \mu_{H_1}^2 > 0$ for all $\mu\in \sosbmh 0 {} X\setminus \{0\}$ by Proposition \ref{M0-char}
 and the just established formula then yields $\snorm \mu_{H}^2 > 0$ for these $\mu$. Consequently, $k$ is characteristic.
 Repeating the argument for $\sosbm X$ yields the last assertion.
\end{proofof}

\begin{proofof}{Lemma \ref{product-form}}
    By the definition of the $\hnorm\cdot$-norm and \eqref{hinner} we have 
 \begin{align*}
  \hnorm {\mu_1\otimes \mu_2}^2 &= \int_{X_1\times X_2} \int_{X_1\times X_2} k_1(x_1,x_1') \cdot k_2(x_2,x_2') \diff\mu_1\!\otimes\! \mu_2(x_1,x_2)\diff \mu_1\!\otimes\! \mu_2(x_1',x_2') \\
&= \int_{X_1\times X_1} \int_{X_2\times X_2} k_1(x_1,x_1') \cdot k_2(x_2,x_2') \diff\mu_1(x_1)\diff \mu_1(x_1')\diff \mu_2(x_2)\diff \mu_2(x_2')\\
 &= \snorm{\mu_1}_{H_1}^2 \cdot \snorm{\mu_2}_{H_2}^2\, ,
 \end{align*}
        where we note that the application of Fubini's theorem is possible, since the kernels are bounded and the measures are finite. Now assume that $k$ is characteristic, but, say, $k_1$ is not strictly integrally positive definite with respect to $\sosbm {X_1}$.
        Then there exists a $\mu_1\in \sosbm{X_1}$ with $\mu_1\neq 0$ but $\snorm{\mu_1}_{H_1}=0$. 
        Moreover, since $\dim \sosbm{X_2}\geq 2$, the decomposition \eqref{m0+p}  gives 
 a $\mu_2\in \sosbmh 0{} {X_2}$ with $\mu_2\neq 0$.
        Let us define $\mu:= \mu_1\otimes \mu_2$. Our construction then yields $\mu\neq 0$
        and $\mu(X_1\times X_2) = \mu_1(X_1) \cdot \mu_2(X_2) = 0$, that is $\mu\in \sosbmh 0{}{X_1\times X_2}$,
 while the already established 
        product formula shows $\hnorm{\mu} = \snorm{\mu_1}_{H_1} \cdot \snorm{\mu_2}_{H_2} = 0$. By Proposition
        \ref{M0-char} we conclude that $k$ is not characteristic.
\end{proofof}

\begin{proofof}{Lemma \ref{plusone-sipd}}
Let $\P_1$ denote the canonical feature map of the kernel $k_1 = \eins_{X\times X}$ and $H_1$ its RKHS.
Moreover, we write  $H+1$ for the RKHS of $k+1$. For later use, we note that 
\eqref{hinner} implies $\langle\P_1(\mu), \P_1(\mu_0)\rangle_{H_1} = 0$ for all 
$\mu\in \sosbm X$ and $\mu_0\in \sosbmh 0 {}X$.

 \atob {iii} {ii} 
 Let us fix  a $P\in \ca M\cap \sosbmh 1 {}X$. Similar to \eqref{m0+p} we then have 
 $\ca M = \R P \oplus \ca M_0$. 
  Indeed, ``$\supset$'' and $\R P \cap \ca M_0=\{0\}$ are trivial and for $\mu \in \ca M$
  it is easy to see that $\mu - \mu(X) P \in \ca M_0$.
 Now, we need to prove $\snorm\mu_{H+1}>0$ for all $\mu = \a P + \mu_0\in \R P \oplus \ca M_0$
 with $\mu\neq 0$.
 By \emph{iii)}  we already know this in the case $\a=0$, and thus we further assume 
 $\a\neq 0$. 
 Then Lemma \ref{char-add} and \eqref{hinner} together with our initial remark
 yield
 \begin{align*}
  \snorm\mu_{H+1}^2 
  &= \snorm\mu_{H}^2 + \langle \P_1(\mu) , \P_1(\mu)\rangle_{H_1}\\
  &= \snorm\mu_{H}^2 + \a^2\snorm P_{H_1}^2 
  + 2\a\langle \P_1(P) , \P_1(\mu_0)\rangle_{H_1}
  + \snorm{\mu_0}_{H_1}^2 \\
  & =  \snorm\mu_{H}^2 + \a^2
 \end{align*}
 and since $\a\neq 0$, we conclude $ \snorm\mu_{H+1}^2 >0$.
 
 \atob {ii} {iii} This is trivial.
 
 \aeqb {iii} i For $\mu_0 \in \ca M_0\subset \sosbmh 0 {}X$,  Lemma \ref{char-add} together with our initial remark shows 
 $\snorm{\mu_0}_{H+1}^2 = \hnorm{\mu_0}^2 + \snorm{\mu_0}_{H_1}^2 = \hnorm{\mu_0}^2$. 
  From this equality the equivalence immediately follows.
\end{proofof}

\begin{proofof}{Theorem \ref{non-equivalence}}
For some fixed
%
 $P\in \sosbmh 1 {} X$ we consider the map 
 \begin{align*}
  \pi : \R P \oplus \sosbmh 0 {} X &\to \sosbm X\\
  \a P + \mu_0 &\mapsto \mu_0\, .
 \end{align*}
  By \eqref{m0+p}, $\pi$ is a linear map $\pi:\sosbm X\to \sosbm X$ with $\pi^2 = \pi$ 
  and $\ran \pi = \sosbmh 0 {} X$,
  and hence a projection onto $\sosbmh 0 {} X$. Moreover, $\sosbmh 0 {} X$ is $\hnorm\cdot$-closed by assumption
  and therefore $\pi$ is   $\hnorm\cdot$-continuous, see \citet[Theorem 3.2.14]{Megginson98}.
  Let us now assume that $\hnorm\cdot$ and $\tvnorm\cdot$ were equivalent on 
  $\sosbmh 0 {} X$.  Our goal is to show that this assumption implies the equivalence of $\hnorm\cdot$ and $\tvnorm\cdot$
	on $\sosbmh {}{}X$.
To this end we fix some sequence $(\mu_n) \in \sosbm X$ and some $\mu \in \sosbm X$
  with $\hnorm{\mu_n-\mu} \to 0$. By \eqref{m0+p} we then find $\a_n, \a\in R$
  with $\mu_n = \a_n P + \pi(\mu_n)$ and $\mu = \a P + \pi(\mu)$.
  Using the $\hnorm\cdot$-continuity of $\pi$, we then find $\hnorm{\pi(\mu_n) - \pi(\mu)} \to 0$, and since
  we assumed that $\hnorm\cdot$ and $\tvnorm\cdot$ are equivalent on 
  $\sosbmh 0 {} X$, we conclude that $\tvnorm{\pi(\mu_n) - \pi(\mu)} \to 0$.
  In addition, we have $\hnorm P>0$ by the injectivity of   the kernel embedding, and therefore
  \begin{displaymath}
   |\a_n-\a| \cdot\hnorm P \leq \hnorm{\mu_n-\mu} + \hnorm{\pi(\mu_n) - \pi(\mu)} \to 0
  \end{displaymath}
  shows that $\a_n\to \a$. Combining these considerations, we find
  \begin{displaymath}
   \tvnorm{\mu_n-\mu} \leq  |\a_n-\a| \cdot\tvnorm P + \tvnorm{\pi(\mu_n) - \pi(\mu)} \to 0\, .
  \end{displaymath}
  Summing up, we have seen that $\hnorm{\mu_n-\mu} \to 0$ implies $\tvnorm{\mu_n-\mu} \to 0$
  for every sequence $(\mu_n) \in \sosbm X$, and consequently, the identity
  map  $\id: (\sosbm X, \hnorm\cdot) \to (\sosbm X, \tvnorm\cdot)$ is continuous.
  This yields $\tvnorm \cdot \leq \snorm\id\cdot\hnorm\cdot$ on $\sosbm X$, and since we also have 
  $\hnorm \cdot \leq \inorm k\cdot\tvnorm\cdot$ we see that 
  $\hnorm\cdot$ and $\tvnorm\cdot$ are indeed  equivalent on 
  $\sosbm X$. This, however,  contradicts  Theorem \ref{norm-equivalence}, and hence our assumption that 
$\hnorm\cdot$ and $\tvnorm\cdot$ were equivalent on 
  $\sosbmh 0 {} X$ is false.
\end{proofof}

\begin{proofof}{Theorem \ref{different-distances}}
 Let us first assume that $\eins_X\in H$. By Lemma \ref{m0-props} we then see that
 the   kernel embedding $\P:\sosbm X\to H$ is injective, and Theorem \ref{non-equivalence} thus
 shows that 
 $\hnorm\cdot$ and $\tvnorm\cdot$ are 
 not equivalent on $\sosbmh 0 {} X$. 
Since $\hnorm\cdot$ is dominated by $\tvnorm\cdot$, we consequently find 
 a sequence $(\mu_n)\subset \sosbmh 0 {} X$
 and some $\d>0$ 
 such that $\hnorm {\mu_n} \to 0$ and  $\inf_{n\geq 1}\tvnorm {\mu_n} \geq \d$.
 Let us consider the Hahn-Jordan decomposition $\mu_n = \mu_n^+ - \mu_n^-$ with $\mu_n^+,\mu_n^-\in \sobm X$.
 Since $\mu_n\in \sosbmh 0 {}X$, we have $\mu_n^+(X) = \mu_n^-(X)$. We define 
 \begin{displaymath}
  Q_n^{(1)} := \frac {\mu_n^+}{\mu_n^+(X)} \qquad \qquad \mbox{ and } \qquad \qquad Q_n^{(2)} := \frac {\mu_n^-}{\mu_n^+(X)}\, .
 \end{displaymath}
  Clearly, this yields $Q_n^{(1)},Q_n^{(2)}\in \sobpm X$ and 
  \begin{displaymath}
   \tvnorm{Q_n^{(1)} -Q_n^{(2)}} = \frac{\mu_n^+(X) + \mu_n^-(X)}{\mu_n^+(X)} = 2\, .
  \end{displaymath}
  Moreover, we have 
    \begin{displaymath}
   \hnorm{Q_n^{(1)} -Q_n^{(2)}} = \frac{\hnorm{\mu_n}}{\mu_n^+(X)} =  \frac{\hnorm{\mu_n}}{\frac{1}{2}\tvnorm{\mu_n}} 
   \leq 2\d^{-1}\hnorm{\mu_n}\, , 
  \end{displaymath}
 and by choosing $n$ sufficiently large we can therefore guarantee  $\hnorm{Q_n^{(1)} -Q_n^{(2)}} \leq \e$.
  
  Let us now consider the case $\eins_X \not\in H$. By Lemma \ref{char-add} we then see that 
  the kernel $\tilde k := k+1$ is characteristic. Let us write $H_1 := \R \eins_X$ for the RKHS of the constant kernel
  $k_1:=\eins_{X\times X}$. Then $x\mapsto (k(x,\cdot), k_1(x,\cdot)) \in H\times H_1$ is a feature map of the kernel
        $k+1$ in the sense of \cite[Definition 4.1]{StCh08} if $H\times H_1$ is equipped with the usual Hilbert space
        norm $\snorm{(w,w_1)} :=\sqrt{ \hnorm w^2 + \snorm{w_1}_{H_1}^2}$ and consequently an application of 
        \cite[Theorem 4.21]{StCh08} shows that the RKHS $\tilde H$ of $k+1$ is given by $\tilde H = H + H_1$. The latter yields
$\eins_X \in \tilde H$.
The already considered case thus shows that for all $\e>0$ 
there exist distributions $Q_1,Q_2\in \sosbmh 1 {} X$ such that 
  $\tvnorm{Q_1-Q_2} = 2$ and $\snorm{Q_1-Q_2}_{\tilde H}\leq \e$. 
Moreover, for $\mu:= Q_1-Q_2 \in   \sosbmh 0 {}X$ 
  we see by Lemma \ref{char-add} and \eqref{hinner}
  that
  \begin{displaymath}
 \e^2 \geq \snorm{Q_1-Q_2}_{\tilde H}^2 =   \snorm\mu_{\tilde H}^2 = \hnorm \mu^2 + \int\int \eins_{X\times X} \diff\mu \diff\mu
   =  \hnorm \mu^2 + \mu(X) \cdot\mu(X) = \hnorm \mu^2 = \hnorm{Q_1-Q_2}^2
  \end{displaymath}
we have shown the assertion
    in the second case, too.
\end{proofof}

\begin {proofof}{Corollary \ref{different-distances-cor}}
   Let us fix some $P\in \sosbmh 1 {} X$,  $\d\in (0,1]$ and $\e\in (0,2\d)$.
	By Theorem \ref{different-distances} there then 
 exist distributions $\tilde Q_1,\tilde Q_2\in \sosbmh 1 {} X$ with 
  $\tvnorm{\tilde Q_1-\tilde Q_2} = 2$ and $\hnorm{\tilde Q_1-\tilde Q_2}\leq \e$.
		We define $Q_i := (1-\d) P + \d \tilde Q_i\in \sosbmh 1 {}X$ for $i=1,2$.
	A simple calculation then shows that
	\begin{align*}
	   \tvnorm{Q_1-Q_2} = \d \tvnorm{\tilde Q_1 - \tilde Q_2} = 2\d\, ,
	\end{align*}
  and analogously, we find 
	$\hnorm{Q_1-Q_2} \leq \d \e \leq \e$. Finally, we have 
$\tvnorm{P- Q_i} = \d \tvnorm{P-\tilde Q_i} \leq 2\d$, and hence we obtain the assertion for $\tilde \d:= 2\d$.
\end {proofof}

\begin{proofof}{Lemma \ref{square-integral-formula}}
Let us first assume that $\mu_1, \mu_2 \in \sobm X$. 
For $x,x'\in X$ we define 
\begin{displaymath}
   \tilde k(x,x') := \sum_{i\in I} \lb_i |e_i(x)| \cdot  |e_i(x')|\, .
\end{displaymath}
Then the Cauchy-Schwarz inequality together with \eqref{kernel_sum_cont} gives
\begin{displaymath}
   \bigl|  \tilde k(x,x') \bigr| 
=  \sum_{i\in I} \lb_i |e_i(x)| \cdot  |e_i(x')|
\leq \sum_{i\in I} \lb_i e_i^2(x) \cdot \sum_{i\in I} \lb_i e_i^2(x')
= k(x,x) \cdot k(x',x')
\end{displaymath}
 for all $x,x'\in X$, and since $k$ is bounded, we conclude that $\tilde k$ is also bounded,
and hence $\tilde k\in \Lx 2 {\mu_1\otimes \mu_2}$. Moreover, for finite $J\subset I$ and $x,x'\in X$ we have 
\begin{displaymath}
   \Bigl| \sum_{i\in J} \lb_i e_i(x)e_i(x')  \Bigr| \leq \tilde k(x,x')\, .
\end{displaymath}
Using Fubini's theorem, \eqref{kernel_sum_cont}, Lebesgue's dominated convergence theorem,
and yet another time  Fubini's theorem,
 we thus find 
\begin{align*}
     \int_X\int_X k(x,x') \, \diff\mu_1(x)\diff\mu_2(x') 
&= \int_{X\times X} \sum_{i\in I} \lb_i e_i(x) e_i(x')\, \diff(\mu_1\!\otimes\! \mu_2)(x,x') \\
&=  \sum_{i\in I} \int_{X\times X} \lb_i e_i(x) e_i(x')\, \diff(\mu_1\!\otimes\! \mu_2)(x,x') \\
& = \sum_{i\in I} \int_X\int_X\lb_i e_i(x) e_i(x') \, \diff\mu_1(x)\diff\mu_2(x') \, .
\end{align*}
From the latter the assertion immediately follows.

Finally, let us assume that $\mu_1,\mu_2 \in \sosbm X$. Using the Hahn decomposition 
 $\mu_1 = \mu_1^+-\mu_1^-$ and $\mu_2=\mu_2^+-\mu_2^-$ 
as well as the fact that the expressions on the left and right of the 
desired equation are linear in the involved measures, we then obtain the assertion by 
the already established case.
\end{proofof}

\begin{proofof}{Theorem \ref{kernel-metric-thm}}
Clearly, the signed measure $\mu:= P-Q$ has the $\nu$-density $h-g$ and we have 
 $\aec{h-g}\in \Ltn \cap \Lx 1 \nu$.
Now we find \eqref{kernel-metric-thm-hxx} by
\eqref{hinner} and Lemma \ref{square-integral-formula}, namely
\begin{align*}
   \g_k^2(P,Q) 
= \snorm{\mu}_H^2
= \int_X\int_X k(x,x') \, \diff\mu(x)\diff\mu(x')
&= \sum_{i\in I} \lb_i \cdot \Bigl(\int_X e_i \, \diff\mu  \Bigr)^2 \\
&= \sum_{i\in I} \lb_i \cdot \Bigl(\int_X e_i \cdot (h-g)  \diff\nu  \Bigr)^2\, .
\end{align*}
%
%

\atob {ii} i 
We split $f$ into $f = f^+ - f^-$ with $f^+\geq 0$ and $f^-\geq 0$.
By our assumption we then know that 
\begin{displaymath}
 c:= \int_X f^+  \diff\nu =  \int_X f^- \diff\nu >0
\end{displaymath}
For $\d>0$, $h\in \D(\nu)$ and $P:= h\diff \nu\in \ca P_2(\nu)$
%
 we define $h_1 := (1+\d c)^{-1}(h+\d f^+)$ and $h_2 := (1+\d c)^{-1}(h+\d f^-)$
and consider the corresponding measures $Q_1 := h_1\diff\nu$ and $Q_2 := h_2\diff\nu$.
The construction immediately ensures $Q_1,Q_2\in \ca P_2(\nu)$. Moreover, we have 
\begin{align*}
 \tvnorm{P-Q_1} = \int_X|h-h_1| \diff\nu
 &= \int_X \Bigl| \frac {h + \d c h} {1+\d c}-  \frac{h+\d f^+}{1+\d c}   \Bigr| \diff\nu \\
 & = \int_X \Bigl| \frac { \d c h  - \d f^+} {1+\d c}  \Bigr| \diff\nu \\
 & \leq \frac \d {1+\d c} \int_X  |  c h| + |f^+|  \diff\nu  \\
 & = \frac {2\d c} {1+\d c} \, ,
\end{align*}
and analogously we find $\tvnorm{P-Q_2}\leq  2\d c/(1+\d c)$. 
In addition, we have 
\begin{align*}
   \tvnorm{Q_1-Q_2} =   \frac 1 {1+\d c}  \int_X   |\d f^+ - \d f^-| \diff\nu
=  \frac \d {1+\d c}  \int_X  |f| \diff\nu
= \frac {2\d c} {1+\d c}
\end{align*}
and by using that $\{2\d c/(1+\d c): \d>0\} = (0,2)$ we obtain the norm conditions for $Q_1$ and $Q_2$.
Finally, \eqref{kernel-metric-thm-hxx} yields
\begin{displaymath}
  \g_k^2(Q_1,Q_2) = \sum_{i\in I} \lb_i  \langle \aec{h_1-h_2}, [e_i]_\sim\rangle_{\Ltn}^2
  = \frac \d {1+\d c}\sum_{i\in I} \lb_i  \langle \aec f, [e_i]_\sim\rangle_{\Ltn}^2 = 0\, ,
\end{displaymath}
which shows \emph{i)} as well as the final assertion.

\atob i {iii} Let $Q_1$ and $Q_2$ be according to \emph{i)} and $\aec{h_1}, \aec{h_2}\in \D(\nu)$ be their
$\nu$-densities. Then $Q_1\neq Q_2$ implies $\aec{h_1}\neq \aec{h_2}$, and \eqref{kernel-metric-thm-hxx}
yields
\begin{displaymath}
   0 =  \g_k^2(P,Q) = \sum_{i\in I} \lb_i  \bigl\langle \aec{h_1-h_2}, [e_i]_\sim \bigr\rangle_{\Ltn}^2  \, .       
\end{displaymath}
 Since $\lb_i> 0$ for all $i\in I$, we then conclude that $\langle \aec{h_1-h_2}, [e_i]_\sim\rangle_{\Ltn}^2=0$ for all 
$i\in I$, which in turn implies \emph{iii)}.

\atob {iii} {ii} We define $f:= h_1-h_2$. Clearly, we have $\aec f\in \Ltn \cap \Lx 1 \nu$ and $f\neq 0$.
Moreover, $h_1,h_2\in \D(\nu)$ gives 
\begin{displaymath}
\int_Xf \diff\nu = \int_X h_1\diff\nu - \int_X h_2\diff\nu = 1-1=0\, .
\end{displaymath}
Finally, the equality 
$\langle\aec{ h_1}, [e_i]_\sim\rangle_{\Ltn} =  \langle \aec{h_2}, [e_i]_\sim\rangle_{\Ltn}$, which holds
 for all $i\in I$,
implies $\langle \aec{f}, [e_i]_\sim\rangle_{\Ltn} = 0 $ for all $i\in I$, and hence $\aec f\in [H]_\sim^\perp$.
\end{proofof}

\begin{proofof}{Corollary \ref{no-uniform}}
   Let $h_{j} := (\eins_X + \a e_j)$ be the $\nu$-density of $Q_{j}$. Using 
 $\inorm{e_j}\leq c_\infty$ and $\a =c_\infty^{-1}$, we then find $h_{j} \geq 0$, and 
since $\aec{e_{i_0}} \perp \aec{e_j}$ we further find 
\begin{displaymath}
   \int_X h_{j}  \diff\nu = 1 + \a \int_X e_j e_{i_0} \diff\nu = 1\, .
\end{displaymath}
Consequently, $Q_{j}$ is   a probability measure. Moreover, we find 
\begin{displaymath}
   \tvnorm{P-Q_j} =  \int_X\bigl|\eins_X - h_{j}\bigr|\diff\nu =  \int_X |\a e_j| \diff\nu \geq c_1  c_\infty^{-1}
\end{displaymath}
and  \eqref{kernel-metric-thm-hxx} implies
\begin{displaymath}
   \snorm{P-Q_j}_H^2 =  \sum_{i\in I} \lb_i  \langle \aec{\eins_X - h_j}, [e_i]_\sim\rangle_{\Ltn}^2= \sum_{i\in I} \lb_i  \langle \aec{- \a e_j}, [e_i]_\sim\rangle_{\Ltn}^2 =  \lb_j  c_\infty^{-2}\, .
\end{displaymath}
This shows the assertion.
\end{proofof}

\begin{proofof}{Corollary \ref{cor:codim}}
In both cases it suffices to show \emph{ii)} of Theorem \ref{kernel-metric-thm}.

\ada i Since $\codim [H]_\sim \geq 2$ there exist linearly independent $\aec {f_1}, \aec {f_2} \in [H]_\sim^\perp$.
If $\int_X f_2 \diff\nu= 0$ then there is nothing left to prove, and if $\int_X f_2 \diff\nu\neq 0$
then a quick calculation shows that 
\begin{displaymath}
   f := f_1  - \frac{\int_X f_1 \diff\nu}{\int_X f_2 \diff\nu} \cdot f_2
\end{displaymath}
is the desired function.

   \ada {ii} Since  $\codim [H]_\sim \geq 1$, there exists an $\aec f \in [H]_\sim^\perp\setminus\{0\}$ 
        and from $\eins_X\in H$ we conclude that 
$\int_X f \diff\nu = \langle \aec f, \aec {\eins_X}\rangle_\Ltn = 0$. 
\end{proofof}

\subsection{Proofs related to Section \ref{sec:lcs}}

\begin{proofof}{Theorem \ref{exist-sipd}}
We fix a $C_0(X)$-kernel $k$  and a finite signed measure
$\mu\in \sosbm X\setminus \sosrm X$. Then $f\mapsto \int_X f\diff\mu$ defines a bounded linear operator 
$C_0(X)\to \R$, and by \eqref{riesz-repres} there thus exists a $\mu^*\in \sosrm X$
such that 
\begin{displaymath}
   \int_X f \diff\mu = \int_X f\diff\mu^*
\end{displaymath}
for all $f\in C_0(X)$. By $H\subset C_0(X)$ and \eqref{alternative-comp} we conclude that 
$\hnorm{\mu-\mu^*}=0$ while our construction ensures $\mu\neq \mu^*$.
\end{proofof}

\begin{proofof}{Theorem \ref{thm-new-char-2}}
%
Using
the already observed identity $I'=\P$ and \citet[Theorem 3.1.17]{Megginson98} 
we see   that $I$ has a dense image if and only if
$\P: \sosrm X\to H$ is injective.
By \eqref{sipd} we then conclude that $k$ is universal if and only if 
$k$ is strictly integrally positive definite with respect to $\sosrm X$. 
%
%
%
%
\end{proofof}

\begin{proofof}{Theorem \ref{exists-univer-char-kern}}
 \aeqb {i} {iii} This has been shown in \citet[Theorem 2]{StHuSc06a}. Moreover, it is well-known that 
 compact metrizable Hausdorff spaces are Polish. The equality $\sosrm X = \sosbm X$ then
 follows from Ulam's theorem, see e.g.~\citet[Lemma 26.2]{Bauer01}.
 
 \atob i {ii} If there exists a universal kernel $k$, then we have already shown that 
 $\sosrm X = \sosbm X$.
 Consequently, $k$
 is   strictly integrally positive definite with respect to $\sosbm X$ by Theorem \ref{thm-new-char-2}, 
 and thus characteristic by Proposition \ref{M0-char}.
 
 \atob {ii} {i} Assume that there exists a characteristic kernel $k$. By Proposition \ref{M0-char} we know 
 that $k$ is 
 strictly integrally positive definite with respect to $\sosbmh 0{} X$.
 Then $k+1$ is a bounded and continuous kernel, 
 which is   strictly integrally positive definite with respect to $\sosbm X$ by Lemma \ref{plusone-sipd}.
 Using $\sosrm X \subset \sosbm X$ and Theorem \ref{thm-new-char-2} we conclude that $k$ is universal.
\end{proofof}

\begin{proofof}{Corollary \ref{product-char}}
   \atob i {ii} Since $(X_1\times X_2, \t_1\otimes \t_2)$ is a compact metrizable space, we have 
                $\sosrm {X_1\times X_2} = \sosbm {X_1\times X_2}$, 
and hence the implication follows from Proposition
                \ref{M0-char} and Theorem \ref{thm-new-char-2}.

        \atob {ii} {iii} Since $(X_1,\t_1)$ and $(X_2,\t_2)$
are assumed to be non-trivial, we find $\dim \sosbm {X_1} \geq 2$ and $\dim \sosbm {X_2} \geq 2$.
Now \emph{iii)}
 follows from Lemma \ref{product-form} and Theorem \ref{thm-new-char-2}.
        
        \atob {iii} i This can be shown by the theorem of Stone-Weierstra\ss, see e.g.~\citet[Lemma A.5]{StThScXXa}
        for details.
\end{proofof}

\begin{proofof}{Theorem \ref{thm-new-char-3}}
Before we begin, we write $E:= \spann\{ e_i: i\in I\}$ and denote the RKHS of $k$ by $H$.

  \aeqb  {ii} {iii}  Via the isomorphism \eqref{riesz-repres} 
 between $C_0(X)'$ and $\sosrm X$
 we easily see that \emph{ii)} is equivalent to 
the statement $\p'_{|E} =0 \implies \p'=0$  for all $\p'\in C_0(X)'$ and by Hahn-Banach's theorem, 
see e.g.~\citet[p.~64, Corollary II.3.13]{DuSc58},
the latter is equivalent to \emph{iii)}.

\atob {iii} {i} The chain of inclusions $E\subset H\subset C_0(X)$ immediately gives the desired implication.

\atob {i} {iii} Clearly, we have $E = \spann\{ \sqrt{\lb_i}e_i: i\in I\}$ and since 
$(\sqrt{\lb_i} e_i)_{i\in I}$ is   an ONB of $H$, we conclude that $\overline E^{H} = H$. 
Let us now fix an $\e>0$ and an $f\in C_0(X)$. Since $k$ is universal, there then exists an $h\in H$
with $\inorm {f-h}\leq \e$ and for this $h$ our initial consideration yields an $e\in E$ with 
$\hnorm{e-h}\leq \e$. Combining both estimates we find
\begin{displaymath}
   \inorm{f-e} 
\leq \inorm {f-h} + \inorm{h-e} 
\leq \e + \inorm k\hnorm{h-e} 
\leq (1+\inorm k)\, \e\, ,
\end{displaymath}
and hence $E$ is dense in $C_0(X)$.
   

To check the last statement, let us assume that $E$ is dense in $C_0(X)$.
Using that $\nu$ is finite we then see that 
$E$ is also dense in $C_0(X)$ with respect to $\snorm\cdot_\Ltn$, and since 
$C_0(X)$ is dense in 
$\sLtn$ by the regularity of $\nu$, see e.g.~\citet[Theorem 29.14]{Bauer01},  
we conclude that $E$ is dense in $\sLtn$. This shows that $[E]_\sim =  \spann\{ [e_i]_\sim: i\in I\}$
is dense in $\Ltn$, and therefore $([e_i]_\sim)_{i\in I}$ is indeed an ONB of $\Ltn$.
\end{proofof}

\begin{proofof}{Corollary \ref{uni-cor}}
 \ada i Let us first assume that  $\lb_i>0$ for all $i\in I$. 
By \citet[Lemma 2.6 and Theorem 2.11]{StSc12a} we then see that 
 $k$ is of the form considered
 in Theorem \ref{thm-new-char-3}, and since  $\spann\{e_i:i\in I\} $ is dense in $C_0(X)$ by our assumption,
 we conclude that $k$ is universal. To show the converse implication we set $I^* := \{i\in I: \lb_i>0\}$
 and assume that $k$ is universal but $I^* \neq I$. By the definition of $I^*$  we have 
 \begin{displaymath}
 k(x,x') = \sum_{i\in I^*} \lb_i e_i(x)e_i(x')
\end{displaymath}
for all $x,x'\in X$, and by \citet[Lemma 2.6 and Theorem 2.11]{StSc12a} we thus conclude that $(e_i)_{i\in I*}$
is the family considered in 
Theorem \ref{thm-new-char-3}. Consequently, this sub-family
$([e_i]_\sim)_{i\in I^*}$ is already an ONB in $\Ltn$, which, however, is impossible for $I^* \neq I$.

\ada{ii} 
Let us first 
 assume that $k$ is characteristic and that there exists an $j\neq i_0$ with $\lb_j=0$. 
We have $[e_j]_\sim \in \Lx 1 \nu$ because $\nu(X) < \infty$. Moreover, 
$[e_j]_\sim \perp [e_i]_\sim$ for all $i\neq j$ implies both $\int_X {e_j} \diff\nu =0$ and 
 $[e_j]_\sim \in [H]_\sim^\perp$,
where in the last step we used \eqref{span_ei}.
Consequently, $k$ cannot be characteristic by Theorem \ref{kernel-metric-thm}.

Conversely, assume that $\lb_i >0$ for all $i\neq i_0$. If $\lb_{i_0}>0$ then $k$ is actually universal by the 
already established part \emph{i)}, and thus characteristic. If $\lb_{i_0}=0$, then the kernel  
\begin{displaymath}
   k(x,x') + 1 = \sum_{i\neq i_0} \lb_i e_i(x)e_i(x') +  e_{i_0}(x)e_{i_0}(x')
\end{displaymath}
is universal by part \emph{i)} and thus $k$ is characteristic by Theorem \ref{thm-new-char-2}, Lemma \ref{plusone-sipd}, 
Proposition \ref{M0-char}, and $\sosrm X = \sosbm X$.

\ada {iii} If $\lb_i >0$ for all $i\in I$, then \emph{i)} shows that $k$ is universal, and by 
 Theorem \ref{thm-new-char-2},
Proposition \ref{M0-char}, and $\sosrm X = \sosbm X$ we conclude that $k$ is characteristic.
To show the converse implication, we assume that 
$k$ is characteristic and there is an $i_{0}\in I$ with $\lb_{i_0}=0$.
Since $([e_i]_\sim)_{i\in I}$ is an ONS in $\Ltn$ and
 $[e_{i_0}]_\sim \not\in \overline{[H]_\sim} = \overline{\spann\{ [e_i]_\sim: \lb_i>0\}}$, see \eqref{span_ei},
we conclude that $[e_{i_0}]_\sim \in [H]_\sim^\perp$. On the other hand, 
$\eins_X \in H$ gives $[\eins_X]_\sim \in [H]_\sim$, and thus we find
\begin{displaymath}
   0 = \langle [e_{i_0}]_\sim, [\eins_X]_\sim\rangle_{\Ltn} = \int_X e_{i_0} \diff\nu\, .
\end{displaymath}
Finally, $[e_{i_0}]_\sim \neq 0$ is obvious and $[e_{i_0}]_\sim \in \Lx 1 \nu$ follows from 
$\nu(X)<\infty$, so that Theorem \ref{kernel-metric-thm} shows that $k$ is not characteristic.
\end{proofof}

\subsection{Proofs related to Section \ref{sec:groups}}

\begin{proofof}{Lemma \ref{real-onb}}
To check that $(e_i^*)_{i\in I}$ is an ONS,
we first observe that the equivalences 
$i=j \Leftrightarrow -i=-j$
and $i=-j \Leftrightarrow -i=j$ imply for $a_i,a_j\in \{-1,1\}$ that
\begin{align*}
 \langle &[ e_i + a_i \bar e_i ]_\sim ,  [e_j + a_j \bar e_j ]_\sim\rangle_{\Lx 2 {G, \C}} \\
& =   \langle [ e_i ]_\sim  + a_i [e_{-i}]_\sim  ,  [e_j ]_\sim + a_j [e_{-j}]_\sim \rangle_{\Lx 2 {G, \C}} \\
& =  \langle [ e_i ]_\sim , [e_{j}]_\sim \rangle +a_j  \langle [ e_i ]_\sim , [e_{-j}]_\sim \rangle +  a_i  \langle [ e_{-i} ]_\sim , [e_{j}]_\sim \rangle + a_ia_j \langle [ e_{-i} ]_\sim , [e_{-j}]_\sim \rangle \\
& = (1+a_ia_j) \d_{i,j} +  (a_i+a_j) \d_{i,-j} \, .
\end{align*}
Let us first consider $i,j\in I_0$. Then we have $i=j$ if and only if $i=-j$, 
that is $\d_{i,j} = \d_{i,-j}$,
and thus we find 
\begin{displaymath}
   \langle [e_i^*]_\sim , [e_j^*]_\sim \rangle_{\Lx 2 G} 
= \frac 1 4 \langle [ e_i +  \bar e_i ]_\sim ,  [e_j + \bar e_j ]_\sim\rangle_{\Lx 2 {G, \C}}
= \d_{i,j}\, .
\end{displaymath}
Similarly, for $i,j\in I_+$ we cannot have $\d_{i,-j}=1$, 
since this would imply $i\in I_-$, and hence we obtain
\begin{displaymath}
   \langle [e_i^*]_\sim , [e_j^*]_\sim \rangle_{\Lx 2 G} 
= \frac 1 2 \langle [ e_i +  \bar e_i ]_\sim ,  [e_j + \bar e_j ]_\sim\rangle_{\Lx 2 {G, \C}}
= \d_{i,j}\, ,
\end{displaymath} 
and for $i,j\in I_-$ we find by an analogous reasoning that
\begin{displaymath}
   \langle [e_i^*]_\sim , [e_j^*]_\sim \rangle_{\Lx 2 G} 
= \frac 1 2 \langle [ e_i -  \bar e_i ]_\sim ,  [e_j - \bar e_j ]_\sim\rangle_{\Lx 2 {G, \C}}
= \d_{i,j}\, .
\end{displaymath} 
For the mixed cases, in which $i$ and $j$ belong to different partition elements $I_0$, $I_+$, or $I_-$
we clearly have $i\neq j$, and hence the above calculation reduces to 
\begin{displaymath}
   \langle [ e_i + a_i \bar e_i ]_\sim ,  [e_j + a_j \bar e_j ]_\sim\rangle_{\Lx 2 {G, \C}}  = (a_i+a_j) \d_{i,-j} \, .
\end{displaymath}
Now, if $i\in I_0$ and $j\not \in I_0$, then we cannot have $i=-j$ and hence we obtain 
$ \langle [e_i^*]_\sim , [e_j^*]_\sim \rangle_{\Lx 2 G} = 0$. For the same reason we find for $i\in I_+$ and $j\in I_-$
with $i\neq -j$ that  $ \langle [e_i^*]_\sim , [e_j^*]_\sim \rangle_{\Lx 2 G} = 0$.
Finally, for $i\in I_+$ and $j\in I_-$ with $i=-j$ we have 
\begin{displaymath}
    \langle [e_i^*]_\sim , [e_j^*]_\sim \rangle_{\Lx 2 G} = 
\frac 1 2 \langle [ e_i +  \bar e_i ]_\sim , \iu\cdot [e_j - \bar e_j ]_\sim\rangle_{\Lx 2 {G, \C}} =
-\iu \cdot(1-1)\cdot \d_{i,-j} = 0\, ,
\end{displaymath}
and therefore we conclude that $([e_i^*]_\sim)_{i\in I}$ is an ONS in $\Lx 2 G$. 

Our next goal is to 
show that $\spann\{ e_i^*: i\in I\}$ is dense in $C(G)$. To this end,
we fix an $f\in C(G)$ and an $\e>0$. Since $f\in C(G,\C)$ and 
$(e_i)_{i\in I}$,
 is dense in $C(G, \C)$, there then exists a finite set $J\subset I$ and $(a_j)_{j\in J}\subset \C$ such that 
\begin{displaymath}
  \sup_{x\in G} \Bigl|  \sum_{j\in J} a_j e_{j}(x) - f(x)\Bigl| < \e\, ,
\end{displaymath}
and from this conclude that 
\begin{align*}
 \sup_{x\in G} \Bigl|  \sum_{j\in J} \re a_j \re e_{j}(x) -  \sum_{j\in J} \im a_j \im e_{j}(x)  - f(x)\Bigl| 
&= \sup_{x\in G} \Bigl| \re \Bigl( \sum_{j\in J} a_j e_{j}(x) \Bigr)- \re f(x)\Bigl| \\
& \leq \sup_{x\in G} \Bigl|  \sum_{j\in J} a_j e_{j}(x) - f(x)\Bigl| \\
&< \e\, .
\end{align*}
In other words, the span of $(\re e_i)_{i\in I} \cup (\im e_i)_{i\in I}$ is dense in $C(G)$.
However, our initial considerations showed that $\re e_i = \re e_{-i}$ and $\im e_i = -\im e_{-i}$ for all 
$i\in I$ as well as $\im e_i =0$ for all $i\in I_0$, and therefore 
$\spann\{ e_i^*: i\in I\}$ is dense in $C(G)$, too.

Finally, $\nu$ is regular, and therefore 
$C(G)$ is dense in $\Lx 2 G$, see e.g.~\citet[Theorem 29.14]{Bauer01}.
Consequently, $\spann\{ e_i^*: i\in I\}$ is dense in $\Lx 2 G$, and therefore $(e_i^*)_{i\in I}$ is an
ONB of $\Lx 2 G$. The estimate $\inorm{e_i^*}\leq \sqrt 2$ follows from $|e_i(x)|= 1$ for all $x\in G$.
\end{proofof}

\begin{proofof}{Lemma \ref{real-bochner}}
Let $\nu$ be the Haar measure of $G$.

\atob {i}{ii}
For a character $e_i\in \hat G$ and $x\in G$ a simple calculation shows
\begin{align*}
   \int_G k(x,x') e_i(x')\diff\nu( {x'}) 
= \int_G \k(-x +x') e_i(x')\diff\nu( {x'} ) 
&= \int_G \k(x') e_i(x+x')\diff\nu ({x'})  \\
&= \lb_i e_i(x)  \, ,
\end{align*}
where $\lb_i :=  \int_G \k(x') e_i(x')\diff\nu ({x'})$. Now, since $k$ is $\R$-valued, the integral operator 
$T_k^\C:\Lx 2 {G,\C} \to \Lx 2 {G,\C}$ is self-adjoint and the previous calculation shows that each character 
$e_i$ gives an eigenvector $[e_i]_\sim\in \Lx 2 {G,\C}$ of $T_k^\C$ with eigenvalue $\lb_i$.  
Using that eigenvalues of self-adjoint operators are real numbers, we then find 
\begin{align*} 
 T_k [\re e_i]_\sim  = T_k^\C [\re e_i]_\sim &= \lb_i [\re e_i]_\sim \\ 
T_k [\im e_i]_\sim= T_k^\C  [\im e_i]_\sim &= \lb_i  [\im e_i]_\sim 
\end{align*}
for all $i\in I$. Now recall that for $i\in I_0$ we have $\im e_i = 0$ and therefore these 
functions
$\im e_i$ 
are \emph{not} eigenvectors of $T_k:\Lx 2 G\to \Lx 2 G$. 
By Lemma \ref{real-onb}
we thus conclude that $([e_i^*]_\sim)_{i\in I}$
is an ONB of eigenvectors of $T_k$ with corresponding 
eigenvalues $(\lb_i)_{i\in I}$, and in particular, there are no further eigenvalues than these.
Moreover, $(\lb_i)_{i\in I}$ is summable since \eqref{int-diag} holds.
Moreover, for $i\in I$ we have $\lb_i = \lb_{-i}$, where we note that for $i\not \in I_0$ the corresponding eigenvalues 
have thus a geometric  multiplicity of at least two.
Since $\nu$ is finite and strictly positive, we thus see that 
$k$ enjoys a Mercer representation \eqref{kernel_sum_cont} for the index set $I^*:= \{i\in I: \lb_i>0\}$
and the sub-family $(e_i^*)_{i\in I^*}$. Moreover, for $x,x'\in G$ this 
representation yields
\begin{align*} \nonumber
   k(x,x')
 &= \sum_{\lb_i > 0} \lb_i e_i^*(x) e_i^*(x') \\ \nonumber
&= \sum_{i\in I_0: \lb_i>0} \lb_i \re e_i(x) \re e_i(x') \\ \nonumber
&\qquad + 2\sum_{i\in I_+: \lb_i>0} \lb_i \bigl( \re e_i(x) \re e_i(x') + \im e_i(x) \im e_i(x') \bigr) \\ \nonumber
& =  \sum_{i\in I_0: \lb_i>0} \lb_i \re e_i(-x'+x) + 2\sum_{i\in I_+: \lb_i>0} \lb_i \re e_i(-x' +x)  \\ \nonumber
& = \sum_{\lb_i> 0} \lb_i \re e_i(-x'+x)\, ,
\end{align*}
where in the second to last step we used \eqref{add-thm} and the last step rests on $\re e_i = \re e_{-i}$.
In addition, $\sup_{i\in I}\inorm{e_i^*}\leq \sqrt 2$  together with the summability of $(\lb_i)_{i\in I}$ 
quickly shows that the series converge both absolutely and uniformly. Finally, the 
continuity of $k$ follows from the uniform convergence in \eqref{Mercer-on-G} and  $e_i^*\in C(G)$
for all $i\in I$.

\atob {ii}{i}
We first observe that Lemma \ref{real-onb} together with 
\citet[Lemma 2.6]{StSc12a} shows that \eqref{Mercer-on-G} does indeed define a kernel $k$, and its translation
invariance is built into the construction. Clearly, $k$ is measurable and
$\sup_{i\in I}\inorm{e_i^*}\leq \sqrt 2$  together with the summability of $(\lb_i)_{i\in I}$ 
shows that $k$ is bounded.
\end{proofof}

\begin{proofof}{Theorem \ref{exist-univ-on-G}}
The equivalences \emph{i)} $\Leftrightarrow$ \emph{iv)} $\Leftrightarrow$ \emph{vi)}
have already been shown in Theorem \ref{exists-univer-char-kern},
and \emph{iii)} $\Rightarrow$ \emph{iv)} is trivial. In addition, if one of the 
conditions \emph{iii)} - \emph{vi)} are satisfied, Theorem \ref{exists-univer-char-kern}
shows that $\sosrm G = \sosbm G$, where in the case of \emph{v)} we additionally need  Lemma \ref{real-bochner}.
Therefore Lemma \ref{plusone-sipd} together with 
Proposition \ref{M0-char} and Theorem \ref{thm-new-char-2} shows that a continuous kernel $k$ is 
characteristic if and only if $k+1$ is universal. This yields the 
equivalences \emph{iv)} $\Leftrightarrow$ \emph{vi)},
and by the last statement in Lemma \ref{real-bochner}, also
 \emph{iii)} $\Leftrightarrow$ \emph{v)}.
%
%
%
It thus suffices to show that 
\emph{ii)} $\Rightarrow$ \emph{iii)} and \emph{i)} $\Rightarrow$ \emph{ii)}.

\atob {ii} {iii} If $\hat G$ is at most countable, so is $I$, and hence there
exists a family $(\lb_i)_{i\in I}$ with $\lb_i>0$ for all $i\in I$ and $\sum_{i\in I}\lb_i<\infty$.
The kernel $k$, which is constructed by \eqref{Mercer-on-G}, is then translation-invariant and continuous, 
and by Corollary \ref{uni-cor} it is also universal.

\atob {i} {ii} 
By \citet[Theorem 3.2.11 and Corollary 3.3.2]{Schurle79} we know that $G$ is completely regular and hence 
\citet[Theorem V.6.6.]{Conway90} shows that $G$ is metrizable if and only if $C(G)$ is separable.
We therefore see that 
%
$C(G)$ is separable.
In addition, since $\nu$ is regular,
\citet[Theorem 29.14]{Bauer01} shows that 
$C(G)$ is dense in $\Lx 2 G$.  
Consequently, $\Lx 2 G$ is separable, 
and since $(e_i^*)_{i\in I}$ is an ONB of $\Lx 2 G$, we conclude that 
  $I$, and thus $\hat G$, is 
at most countable.
\end{proofof}

\begin{proofof}{Corollary \ref{char-on-G-char}}
 Since the Haar measure $\nu$ on $G$ is a finite, regular, and strictly positive Borel measure, 
 we see by Lemma \ref{real-onb} that all assumptions of Corollary \ref{uni-cor} are satisfied.
 Moreover, we have $e^*_0 = \eins_G$. Now the assertions follow from Corollary \ref{uni-cor}.
\end{proofof}

\subsection{Proofs related to Section \ref{sec:Sd}}

\begin{proofof}{Proof of Theorem \ref{thm:32}}
We first consider the 
case $\psi \in \Psi_{d+2}$. If the kernel $k$ on $\bbS^d$ induced by $\psi$ is strictly positive definite then Lemma \ref{prop:1} implies that $b_{n,d}>0$ for all $n \ge 0$. By Theorem \ref{thm:Schoenberg}, $k$ is then universal and characteristic. Conversely, if $k$ is universal or characteristic then $b_{n,d} > 0$ for all $n \ge 1$ by Theorem \ref{thm:Schoenberg}, thus it is strictly positive definite by Remark \ref{rem:spd2cb}.

In the case $\psi \in \Psi_{d+1}^+$, we have $\psi \in \Psi_{d}^+$, and hence it suffices to show that
$k$ is universal and characteristic. This, however, follows from Lemma \ref{prop:34} and Theorem \ref{thm:Schoenberg}.
\end{proofof}

\begin{proofof}{Theorem \ref{thm:33}}
Suppose that the kernel $k$ is induced by $\psi\in \Psi_{\infty}^+$. Then $b_{n,d} >0$ for all $n \in \mathbb{N}_0$
by Lemma \ref{prop:34} and hence the kernel is universal and characteristic by Theorem \ref{thm:Schoenberg}. Suppose now that $k$ is characteristic. By Proposition \ref{prop:421} we obtain that $\psi \in \Psi_{\infty}^+$. 
\end{proofof}

\begin{proofof}{Lemma \ref{prop:1}}
It is clear from the results of \citep{ChenMenegattoETAL2003} that  $b_{n,d} > 0$ for all $n \in \mathbb{N}_0$ is a sufficient condition for $\psi$ being strictly positive definite. Suppose now that $\psi \in \Psi_{d+2} \cap \Psi_{d}^+$. \citet[Corollary 4]{Gneiting2013} implies that if $b_{2k+2,d} > 0$ ($b_{2k+1,d}> 0$) for some $k$, then $b_{2k' + 2,d} > 0$ ($b_{2k'+1,d}> 0$) for all $k' \le k$. This yields the claim by Remark \ref{rem:spd2cb}.
\end{proofof}

\begin{proofof}{Lemma \ref{prop:34}}
This is shown in \citet[Proof of Corollary 1(b)]{Gneiting2013}. 
\end{proofof}

\begin{proofof}{Proposition \ref{prop:421}}
Assume that $\psi$ is not strictly positive definite, or, by Remark \ref{rem:spd2cb}, does not satisfy condition $b$. We will show that it cannot be characteristic.

First, we construct a special class of probability densities $p$ on $\bbS^d$ such that we explicitly know the integrals 
\begin{equation}\label{eq:ckjp}
c_{k,j}(p) := \int_{\bbS^d} e_{k,j}(y) p(y)\diff \sigma(y) 
\end{equation}
with respect to the basis of spherical harmonics. Here, $\sigma$ is the surface area measure on $\mathbb{S}^d$ normalized such that $\int_{\mathbb{S}^d}\diff \sigma = 1$. 
Fix $v_0 \in \mathbb{S}^d$ and $a \in [-1,1]\backslash \{0\}$. We have 
\begin{align*}
 |C_n^{(d-1)/2}(x)|\le C_n^{(d-1)/2}(1)\, , \qquad \qquad   x \in [-1,1],
\end{align*}
see   (\citetalias[18.14.4]{dlmf} for $d \ge 2$), and
therefore, for all $n \ge 1$, $x \in \bbS^d$, we have
\begin{equation}\label{eq:pna}
p_{n,a}(x) := 1 + a\frac{C_n^{(d-1)/2}(\langle v_0,x\rangle)}{C_n^{(d-1)/2}(1)} \ge 0,
\end{equation}
and $\int_{\bbS^d} p_{n,a}(x)\diff \sigma(x) = 1$,
where for the last equality we used that $C_{n}^{(d-1)/2}(\langle v_0,\cdot\rangle)$ is a spherical
harmonics of degree $n$ and thus orthogonal to $e_{0,0} = \eins_{\bbS^d}$.
Consequently,  $p_{n,a}$ is a probability density function on $\bbS^d$ 
with respect to the surface area measure $\sigma$. Note that $p_{n,a}$ and $p_{n',a}$ induce different probability measures on $\mathbb{S}^d$ for $n \not=n'$. We obtain 
\begin{equation}\label{eq:ckjpna}
c_{k,j}(p_{n,a}) = \delta_{k,0} + \delta_{k,n}\frac{a}{N(d,k)}e_{k,j}(v_0)
\end{equation}
using \eqref{eq:legendre}, and 
the Funk-Hecke Theorem \citep[Theorem 3.4.1]{Groemer1996} yields that
\begin{equation*}
\int_{\bbS^d} \langle x,y \rangle^n e_{k,j}(y) \diff \sigma(y) = \lambda_k^n e_{k,j}(x), \quad x \in \bbS^d, j = 1,\dots,N(d,k),
\end{equation*}
where
\[
\lambda_k^n = \frac{\Gamma((d+1)/2)}{\sqrt{\pi}\Gamma(d/2)}C_k^{(d-1)/2}(1)^{-1}\int_{-1}^1 t^n C_k^{(d-1)/2}(t)(1-t^2)^{(d-2)/2}dt.
\]
Since the family $(e_{n,j})_{n\in \N_o, j=1,\dots,N(d,k)}$ is an ONB of $\Lx 2 {\bbS^d}$,
we obtain the following  Mercer representation
of the bounded and continuous kernel $(x,y)\mapsto \langle x,y\rangle^n$ on $\bbS^d$ 
\begin{equation}\label{eq:xyton}
\langle x,y\rangle^n = \sum_{k=0}^{\infty} \sum_{j=1}^{N(d,k)} \lambda_k^n e_{k,j}(x)e_{k,j}(y)\, ,
\end{equation}
where for each $x$ the convergence is uniform in $y$ by \citet[Corollary 3.5]{StSc12a}.
Using that $C_k^{(d-1)/2}$ is even for $k$ even and odd for $k$ odd, one obtains that $\lambda_k^n=0$ if $k-n$ is odd. The $C_k^{(d-1)/2}$ are orthogonal with respect to the weight function $(1-t^2)^{(d-2)/2}$ \citetalias[18.3.1]{dlmf}, therefore $\lambda_k^n = 0$ for $k > n$. Finally, the formula \citetalias[18.17.37]{dlmf} for the Mellin transform yields that 
\[
\lambda_k^n = \frac{\pi 2^{d-n-1}\Gamma(k+d-1)\Gamma(n+1)}{k! \Gamma(\frac{d-1}{2})\Gamma(\frac{k+d+n}{2})\Gamma(\frac{n-k + 2}{2})} > 0 \quad k \le n, \; k-n \; \text{even}.
\]
For a probability density $p$ on $\bbS^d$, we have by \eqref{eq:Schoenberg} and \eqref{eq:xyton} for $x \in \bbS^d$,
\begin{align*}
\int_{\bbS^d} k(x,y) p(y) \diff\sigma(y) &= \sum_{n=0}^{\infty} \sum_{k=0}^n \sum_{j=1}^{N(d,k)} b_n\lambda_k^n e_{k,j}(x) c_{k,j}(p)\\
&= \sum_{k=0}^{\infty} z_k \sum_{j=1}^{N(d,k)} c_{k,j}(p) e_{k,j}(x),
\end{align*}
where $z_k = \sum_{n=k}^{\infty} b_n \lambda_k^n$ and $c_{k,j}(p)$ is defined at \eqref{eq:ckjp}. If $b_n = 0$ for all even $n \ge n_0$ then $z_k = 0$ for all even $k \ge n_0$. If $b_n =0$ for all odd $n \ge n_0$ then $z_k = 0$ for all odd $k \ge n_0$. 

Let us start with the case that $b_n = 0$ for all even $n \ge n_0$, i.e.~$z_k = 0$ for all even $k\ge n_0$. For all $m \in \mathbb{N}_0$, we have $c_{k,j}(p_{2m,a}) = 0$ for $k$ odd by \eqref{eq:ckjpna}, where $p_{n,a}$ is defined at \eqref{eq:pna}. Hence, for $2m \ge n_0$ and $x\in \bbS^d$, we obtain 
\begin{align*}
\int_{\bbS^d} &k(x,y) p_{2m,a}(y) \diff\sigma(y) \\&= \sum_{k=0,k\text{ even}}^{n_0} z_k  \sum_{j=1}^{N(d,k)}\left(\delta_{k,0} + \delta_{k,2m}\frac{a}{N(d,k)}e_{k,j}(v_0)\right) e_{k,j}(x) = z_0,
\end{align*}
which shows that the kernel mean embedding maps all these densities to the constant function with value $z_0$.
Consequently, $k$ is not characteristic.

Suppose now that $b_n = 0$ for all odd $n \ge n_0$, i.e.~$z_k = 0$ for all odd $k \ge n_0$. For all $m \in \mathbb{N}_0$, we have $c_{k,j}(p_{2m+1,a}) = \delta_{k,0} $ for $k$ even by \eqref{eq:ckjpna}. Hence, for $2m + 1 \ge n_0$
and $x\in \bbS^d$, we obtain
\begin{align*}
\int_{\bbS^d} &k(x,y) p_{2m+1,a}(y) \diff\sigma(y) \\&= z_0 + \sum_{k=1,k\text{ odd}}^{n_0} z_k  \sum_{j=1}^{N(d,k)} \delta_{k,2m+1}\frac{a}{N(d,k)}e_{k,j}(v_0) e_{k,j}(x)= z_0,
\end{align*}
which again shows that the kernel mean embedding maps all these densities to the constant function with value $z_0$.
Consequently, $k$ is not characteristic.
\end{proofof}

\small{

}


\begin{thebibliography}{42}
\providecommand{\natexlab}[1]{#1}
\providecommand{\url}[1]{\texttt{#1}}
\expandafter\ifx\csname urlstyle\endcsname\relax
  \providecommand{\doi}[1]{doi: #1}\else
  \providecommand{\doi}{doi: \begingroup \urlstyle{rm}\Url}\fi

\bibitem[Bauer(2001)]{Bauer01}
H.~Bauer.
\newblock \emph{Measure and Integration Theory}.
\newblock De Gruyter, Berlin, 2001.

\bibitem[Berg et~al.(1984)Berg, Christensen, and Ressel]{BergChristensETAL1984}
C.~Berg, J.~P.~R. Christensen, and P.~Ressel.
\newblock \emph{Harmonic Analysis on Semigroups}.
\newblock Springer, New York, 1984.

\bibitem[Bochner(1941)]{Bochner1941}
S.~Bochner.
\newblock Hilbert distances and positive definite functions.
\newblock \emph{Ann. of Math.}, 42:\penalty0 647--656, 1941.

\bibitem[Chen et~al.(2003)Chen, Menegatto, and Sun]{ChenMenegattoETAL2003}
D.~Chen, V.~A. Menegatto, and X.~Sun.
\newblock A necessary and sufficient condition for strictly positive definite
  functions on spheres.
\newblock \emph{Proc. Amer. Mat. Soc.}, 131:\penalty0 2733--2740, 2003.

\bibitem[Chen et~al.(2016)Chen, Wang, and Zhang]{ChWaZh16a}
W.~Chen, B.~Wang, and H.~Zhang.
\newblock Universalities of reproducing kernels revisited.
\newblock \emph{Appl. Anal.}, 95:\penalty0 1776--1791, 2016.

\bibitem[Conway(1990)]{Conway90}
J.~B. Conway.
\newblock \emph{A Course in Functional Analysis}.
\newblock Springer-Verlag, New York, 2nd edition, 1990.

\bibitem[Diestel and Uhl(1977)]{DiUh77}
J.~Diestel and J.~J. Uhl.
\newblock \emph{Vector Measures}.
\newblock American Mathematical Society, Providence, 1977.

\bibitem[Diestel et~al.(1995)Diestel, Jarchow, and Tonge]{DiJaTo95}
J.~Diestel, H.~Jarchow, and A.~Tonge.
\newblock \emph{Absolutely Summing Operators}.
\newblock Cambridge University Press, Cambridge, 1995.

\bibitem[{Digital Library of Mathematical Functions (DLMF)}(2011)]{dlmf}
{Digital Library of Mathematical Functions (DLMF)}.
\newblock \emph{Release date 2011-08-29.}
\newblock National Institute of Standards and Technology from
  \url{http://dlmf.nist.gov/}, 2011.

\bibitem[Dunford and Schwartz(1958)]{DuSc58}
N.~Dunford and J.~T. Schwartz.
\newblock \emph{Linear Operators, Part I: General Theory}.
\newblock Interscience Publishers, Inc., New York, 1958.

\bibitem[Folland(1995)]{Folland95}
G.~B. Folland.
\newblock \emph{A course in Abstract Harmonic Analysis}.
\newblock CRC Press, Boca Raton, FL, 1995.

\bibitem[Fukumizu et~al.(2008)Fukumizu, Gretton, Sun, and
  Sch{\"o}lkopf]{FuGrSuSc08a}
K.~Fukumizu, A.~Gretton, X.~Sun, and B.~Sch{\"o}lkopf.
\newblock Kernel measures of conditional dependence.
\newblock In J.~C. Platt, D.~Koller, Y.~Singer, and S.~T. Roweis, editors,
  \emph{Advances in Neural Information Processing Systems 20}, pages 489--496.
  2008.

\bibitem[Fukumizu et~al.(2009)Fukumizu, Sriperumbudur, Gretton, and
  Sch\"{o}lkopf]{FuSrGrSc09a}
K.~Fukumizu, B.~Sriperumbudur, A.~Gretton, and B.~Sch\"{o}lkopf.
\newblock Characteristic kernels on groups and semigroups.
\newblock In D.~Koller, D.~Schuurmans, Y.~Bengio, and L.~Bottou, editors,
  \emph{Advances in Neural Information Processing Systems 21}, pages 473--480.
  Curran Associates, Inc., 2009.

\bibitem[Gneiting(2013)]{Gneiting2013}
T.~Gneiting.
\newblock Strictly and non-strictly positive definite functions on spheres.
\newblock \emph{Bernoulli}, 19:\penalty0 1327--1349, 2013.

\bibitem[Gneiting and Katzfuss(2014)]{GneitingKatzfuss2014}
T.~Gneiting and M.~Katzfuss.
\newblock Probabilistic forecasting.
\newblock \emph{Ann.\ Rev.\ Stat.\ Appl.}, 1:\penalty0 125--151, 2014.

\bibitem[Gneiting and Raftery(2007)]{GneitingRaftery2007}
T.~Gneiting and A.~E. Raftery.
\newblock Strictly proper scoring rules, prediction, and estimation.
\newblock \emph{J. Amer. Statist. Assoc.}, 102:\penalty0 359--378, 2007.

\bibitem[Gretton et~al.(2007)Gretton, Borgwardt, Rasch, Sch\"{o}lkopf, and
  Smola]{GrBoRaScSm07a}
A.~Gretton, K.M. Borgwardt, M.~Rasch, B.~Sch\"{o}lkopf, and A.J. Smola.
\newblock A kernel method for the two-sample-problem.
\newblock In P.~B. Sch\"{o}lkopf, J.~C. Platt, and T.~Hoffman, editors,
  \emph{Advances in Neural Information Processing Systems 19}, pages 513--520.
  MIT Press, 2007.

\bibitem[Groemer(1996)]{Groemer1996}
H.~Groemer.
\newblock \emph{Geometric applications of {F}ourier series and sperical
  harmonics}.
\newblock Cambridge University Press, New York, 1996.

\bibitem[Hewitt and Ross(1963)]{HeRo63}
E.~Hewitt and K.A. Ross.
\newblock \emph{Abstract Harmonic Analysis. {V}ol. {I}}.
\newblock Springer-Verlag, Berlin, 1963.

\bibitem[Hewitt and Ross(1970)]{HeRo70}
E.~Hewitt and K.A. Ross.
\newblock \emph{Abstract Harmonic Analysis. {V}ol. {II}}.
\newblock Springer-Verlag, Berlin, 1970.

\bibitem[Hewitt and Stromberg(1965)]{HeSt65}
E.~Hewitt and K.~Stromberg.
\newblock \emph{Real and Abstract Analysis. {A} Modern Treatment of the Theory
  of Functions of a Real Variable}.
\newblock Springer-Verlag, New York, 1965.

\bibitem[Hofmann and Morris(2013)]{HoMo13}
K.H. Hofmann and S.A. Morris.
\newblock \emph{The Structure of Compact Groups: A Primer for the Student -- A
  Handbook for the Expert}.
\newblock De Gruyter, Berlin, 3rd edition, 2013.

\bibitem[Megginson(1998)]{Megginson98}
R.~E. Megginson.
\newblock \emph{An Introduction to {B}anach Space Theory}.
\newblock Springer-Verlag, New York, 1998.

\bibitem[Menegatto(1992)]{Menegatto1992}
V.~A. Menegatto.
\newblock \emph{Interpolation on spherical spaces}.
\newblock PhD thesis, University of Texas at Austin, 1992.

\bibitem[Menegatto(1994)]{Menegatto1994}
V.~A. Menegatto.
\newblock Strictly positive definite kernels on the {H}ilbert sphere.
\newblock \emph{Appl. Anal.}, 55:\penalty0 91--101, 1994.

\bibitem[Menegatto(1995)]{Menegatto1995}
V.~A. Menegatto.
\newblock Strictly positive definite kernels on the circle.
\newblock \emph{Rocky Mountain J. Math.}, 25:\penalty0 1149--1163, 1995.

\bibitem[Menegatto et~al.(2006)Menegatto, Oliveira, and
  Peron]{MenegattoOliveiraETAL2006}
V.~A. Menegatto, C.~P. Oliveira, and A.~P. Peron.
\newblock Strictly positive definite kernels on subsets of the complex plane.
\newblock \emph{Comput. Math. Appl.}, 51:\penalty0 1233--1250, 2006.

\bibitem[Micchelli et~al.(2006)Micchelli, Xu, and Zhang]{MiXuZh06a}
C.~A. Micchelli, Y.~Xu, and H.~Zhang.
\newblock Universal kernels.
\newblock \emph{J. Mach. Learn. Res.}, 7:\penalty0 2651--2667, 2006.

\bibitem[Morris(1979)]{Morris79a}
S.A. Morris.
\newblock Duality and structure of locally compact {A}belian groups{$\ldots
  $}for the layman.
\newblock \emph{Mathematical Chronicle}, 8:\penalty0 39--56, 1979.

\bibitem[Muandet et~al.(2017)Muandet, Fukumizu, Sriperumbudur, and
  Sch{\"o}lkopf]{MuFuSrSc17a}
K.~Muandet, K.~Fukumizu, B.~Sriperumbudur, and B.~Sch{\"o}lkopf.
\newblock Kernel mean embedding of distributions: A review and beyond.
\newblock \emph{Foundations and Trends in Machine Learning}, 10:\penalty0
  1--141, 2017.

\bibitem[Pillai et~al.(2007)Pillai, Wu, Liang, Mukherjee, and
  Wolpert]{PiWuLiMuWo07a}
N.S. Pillai, Q.~Wu, F.~Liang, S.~Mukherjee, and R.L. Wolpert.
\newblock Characterizing the function space for {B}ayesian kernel models.
\newblock \emph{J. Mach. Learn. Res.}, 8:\penalty0 1769--1797, 2007.

\bibitem[Schoenberg(1942)]{Schoenberg1942}
I.~J. Schoenberg.
\newblock Positive definite functions on spheres.
\newblock \emph{Duke Math. J.}, 9:\penalty0 96--108, 1942.

\bibitem[Schurle(1979)]{Schurle79}
A.W. Schurle.
\newblock \emph{Topics in Topology}.
\newblock Elsevier North Holland, 1979.

\bibitem[Simon-Gabriel and Sch{\"o}lkopf(2016)]{SGScXXa}
C.-J. Simon-Gabriel and B.~Sch{\"o}lkopf.
\newblock Kernel distribution embeddings: Universal kernels, characteristic
  kernels and kernel metrics on distributions.
\newblock Technical report, 2016.
\newblock URL \url{https://arxiv.org/abs/1604.05251}.

\bibitem[Sriperumbudur et~al.(2010{\natexlab{a}})Sriperumbudur, Gretton,
  Fukumizu, Sch\"olkopf, and Lanckriet]{SriperumbGrettonETAL2010}
B.~K. Sriperumbudur, A.~Gretton, K.~Fukumizu, B.~Sch\"olkopf, and G.~R.~G.
  Lanckriet.
\newblock Hilbert space embeddings and metrics on probability measures.
\newblock \emph{J. Mach. Learn. Res.}, 11:\penalty0 1517--1561,
  2010{\natexlab{a}}.

\bibitem[Sriperumbudur et~al.(2011)Sriperumbudur, Fukumizu, and
  Lanckriet]{SriperumbFukumizuETAL2011}
B.~K. Sriperumbudur, K.~Fukumizu, and G.~R.~G. Lanckriet.
\newblock Universality, characteristic kernels and {RKHS} embedding of
  measures.
\newblock \emph{J. Mach. Learn. Res.}, 12:\penalty0 2389--2410, 2011.

\bibitem[Sriperumbudur et~al.(2010{\natexlab{b}})Sriperumbudur, Fukumizu, and
  Lanckriet]{SrFuLa10a}
B.K. Sriperumbudur, K.~Fukumizu, and G.~Lanckriet.
\newblock On the relation between universality, characteristic kernels and rkhs
  embedding of measures.
\newblock volume~9, pages 773--780, 2010{\natexlab{b}}.

\bibitem[Steinwart(2001)]{Steinwart01a}
I.~Steinwart.
\newblock On the influence of the kernel on the consistency of support vector
  machines.
\newblock \emph{J. Mach. Learn. Res.}, 2:\penalty0 67--93, 2001.

\bibitem[Steinwart and Christmann(2008)]{StCh08}
I.~Steinwart and A.~Christmann.
\newblock \emph{Support Vector Machines}.
\newblock Springer, New York, 2008.

\bibitem[Steinwart and Scovel(2012)]{StSc12a}
I.~Steinwart and C.~Scovel.
\newblock {M}ercer's theorem on general domains: on the interaction between
  measures, kernels, and {RKHS}s.
\newblock \emph{Constr. Approx.}, 35:\penalty0 363--417, 2012.

\bibitem[Steinwart et~al.(2006)Steinwart, Hush, and Scovel]{StHuSc06a}
I.~Steinwart, D.~Hush, and C.~Scovel.
\newblock Function classes that approximate the {B}ayes risk.
\newblock In G.~Lugosi and H.~U. Simon, editors, \emph{Proceedings of the 19th
  Annual Conference on Learning Theory}, pages 79--93, New York, 2006.
  Springer.

\bibitem[Steinwart et~al.(2016)Steinwart, Thomann, and Schmid]{StThScXXa}
I.~Steinwart, P.~Thomann, and N.~Schmid.
\newblock Learning with hierarchical {G}aussian kernels.
\newblock Technical report, Fakult{\"a}t f{\"u}r Mathematik und Physik,
  Universit{\"at} Stuttgart, 2016.
\newblock URL \url{http://arxiv.org/abs/1612.00824}.

\end{thebibliography}
\end{document}